\documentclass{amsart}
\usepackage[english]{babel}
\usepackage{enumerate}
\usepackage{amsmath,amssymb,amsthm}
\usepackage [cmtip,arrow]{xy}
\xyoption{all}
\usepackage {pb-diagram,pb-xy}

\ifx\pdfpageheight\undefined
   \usepackage[dvips,colorlinks=true,linkcolor=blue,citecolor=red,%
      urlcolor=green]{hyperref}
   \usepackage[dvips]{graphicx}
   \makeatletter
   \edef\Gin@extensions{\Gin@extensions,.mps}
   \makeatother
\else
 \usepackage[pdftex]{graphicx}
   \usepackage[bookmarksopen=false,pdftex=true,breaklinks=true,%
      backref=page,pagebackref=true,plainpages=false,%
      hyperindex=true,pdfstartview=FitH,colorlinks=true,%
      pdfpagelabels=true,colorlinks=true,linkcolor=blue,%
      citecolor=red,urlcolor=green,hypertexnames=false%
      ]%
   {hyperref}
\fi
\usepackage{bm}
\usepackage{tabularx}
\usepackage{color, colortbl}
\usepackage{url}

\catcode`\<=\active \def<{
\fontencoding{T1}\selectfont\symbol{60}\fontencoding{\encodingdefault}}
\catcode`\>=\active \def>{
\fontencoding{T1}\selectfont\symbol{62}\fontencoding{\encodingdefault}}
\catcode`\|=\active \def|{
\fontencoding{T1}\selectfont\symbol{124}\fontencoding{\encodingdefault}}

\newcommand{\tmmathbf}[1]{\ensuremath{\boldsymbol{#1}}}
\newcommand{\tmop}[1]{\ensuremath{\operatorname{#1}}}

\definecolor{LightCyan}{rgb}{0.88,1,1}

\newtheorem{theorem}{Theorem}
\newtheorem{lemma}{Lemma}[section]
\newtheorem{theorem-corollary}[theorem]{Corollary}
\newtheorem{corollary}[lemma]{Corollary}
\newtheorem{proposition}[lemma]{Proposition}

\theoremstyle{definition}
\newtheorem{definition}[lemma]{Definition}
\newtheorem{example}[lemma]{Example}

\newtheorem{notation}[lemma]{Notation}

\newtheorem{remark}[lemma]{Remark}

\newcommand{\ra}{\rangle}
\newcommand{\la}{\langle}

\newcommand{\R}{\mathrm{R}}
\newcommand{\C}{\mathrm{C}}

\newcommand{\PP}{\mathbb{P}}
\newcommand{\ZZ}{\mathrm{Zer}}

\newcommand{\eps}{\varepsilon}

\DeclareMathOperator{\Ext}{Ext}

\newcommand{\Bas}{\mathrm{Bas}}
\newcommand{\hide}[1]{}

\begin{document}
\title[On a real analogue of Bezout inequality] {On a real analogue of Bezout inequality and  the number of connected
components of sign conditions}

\author[Barone]{Sal Barone}
\address{School of Mathematics\\
Georgia Institute of Technology, Atlanta, GA 30332}
\email{sbarone@math.gatech.edu}
\urladdr{\url{http://people.math.gatech.edu/~sbarone7/}}

\author[Basu]{Saugata Basu}
\address{Department of Mathematics\\
Purdue University, West Lafayette, IN 47907}
\email{sbasu@math.purdue.edu}
\urladdr{\url{http://www.math.purdue.edu/~sbasu}}

\thanks{The second author was partially supported by NSF grants
  CCF-0915954,  CCF-1319080 and DMS-1161629.  }

{\maketitle}

\begin{abstract}
  Let $\R$ be a real closed field and $Q_1,
  \ldots, Q_{\ell} \in \R [X_1, \ldots, X_k]$
  such that for each $i, 1 \leq i \leq \ell$, $\deg (Q_i) \leq d_i$. For $1
  \leq i \leq \ell$, denote by $\mathcal{Q}_i = \{Q_1, \ldots, Q_i \}$, $V_i$
  the real variety defined by $\mathcal{Q}_i$, and $k_i$ an upper bound on the
  real dimension of $V_i$ (by convention $V_0 =
  \R^k$ and $k_0 = k$). Suppose also that
  \[ 2 \leq d_1 \leq d_2 \leq \frac{1}{k + 1} d_3 \leq \frac{1}{(k + 1)^2}
     d_4 \leq \cdots \leq \frac{1}{(k + 1)^{\ell - 3}} d_{\ell - 1} \leq
     \frac{1}{(k + 1)^{\ell - 2}} d_{\ell}, \]
  and that $\ell \leq k$. We prove that the number of semi-algebraically
  connected components of $V_{\ell}$ is bounded by
  \[ O (k)^{2 k} \left( \prod_{1 \leq j < \ell} d_j^{k_{j - 1} - k_j}
     \right) d_{\ell}^{k_{\ell - 1}} . \]
  This bound can be seen as a weak extension of the classical Bezout
  inequality (which holds only over algebraically closed fields and is
  false over real closed fields) to varieties defined over real
  closed fields.
  
  Additionally, if $\mathcal{P} \subset \R [X_1,
  \ldots, X_k] $ is a finite family of polynomials with $\deg (P) \leq d$ for
  all $P \in \mathcal{P}$, $\ensuremath{\operatorname{card}}\mathcal{P}= s$,
  and $d_{\ell} \leq \frac{1}{k + 1} d$, we prove that the number of
  semi-algebraically connected components of the realizations of all
  realizable sign conditions of the family $\mathcal{P}$ restricted to
  $V_{\ell}$ is bounded by
  \[ O (k)^{2 k} (s d)^{k_{\ell}}  \left( \prod_{1 \leq j \leq \ell} d_j^{k_{j
     - 1} - k_j} \right) . \]
These results have found applications in discrete geometry, for proving incidence bounds
\cite{Basu-Sombra2014},  as well as in effcient range-searching \cite{Matousek-Safernova2014}.
\end{abstract}

\tableofcontents

\section{Introduction }
\label{sec:intro}

\subsection{History and motivation}
Let $\R$ be a fixed real closed field, and we
denote by $\C$ the algebraic closure of
$\R$. Bounds on the number of semi-algebraically
connected components, and in fact on all the Betti numbers of real algebraic
varieties and of semi-algebraic subsets of
$\R^k$ in terms of the number and the degrees of
the polynomials used to define them is a well studied problem in quantitative
real algebraic geometry. The classical bounds, going back to the work of
 Ole{\u\i}nik and Petrovski{\u\i} {\cite{OP}}, Thom {\cite{T}} and Milnor
{\cite{Milnor2}}, bounded the sum of the Betti numbers of real algebraic
varieties, as well as those of basic closed semi-algebraic sets. These and
related bounds (see below) are extremely important in real algebraic geometry
{\cite{BPRbook2}}, but have also been used extensively in other areas such
as combinatorics {\cite{Alon}}, discrete and computational geometry
{\cite{GPW96}}, and theoretical computer science {\cite{montana-pardo-1993}} \
(the cited references are not by any means exhaustive but only given for
illustrative purposes -- we refer the reader to {\cite{BPR10}} for a more
extensive survey).

An important application of the bounds mentioned above is in bounding the
number of semi-algebraically connected components of the realizations of
various sign conditions of a family of polynomials in
$\R^k$ or more generally sign conditions
restricted to a given real sub-variety of
$\R^k$. In order to state these results more
precisely, we introduce some notation.

\begin{notation}
  For $\mathcal{P} \subset \R [X_1, \ldots,
  X_k]$ a finite family of polynomials, a {{\em sign condition\/}} $\sigma$ on
  $\mathcal{P}$ is an element of $\{ 0, 1, - 1 \}^{\mathcal{P}}$. The
  \emph{realization $\ensuremath{\operatorname{Reali}} (\sigma, V)$ of the sign
  condition $\sigma$ on a semi-algebraic set $V\subset \R^k$} is the semi-algebraic
  set defined by
  \begin{eqnarray*}
    \ensuremath{\operatorname{Reali}} (\sigma, V) & = & \{ x \in V \mid
    \ensuremath{\operatorname{sign}} (P) = \sigma (P) , P \in \mathcal{P} \} .
  \end{eqnarray*}
\end{notation}

\begin{notation}
  For any finite family of polynomials $\mathcal{Q} \subset
  \R [X_1, \ldots, X_k]$ we will denote by
  $\ensuremath{\operatorname{Zer}} (\mathcal{Q},
  \R^k)$ the set of real zeros of $\mathcal{Q}$
  in $\R^k$. If $\mathcal{Q}= \{ Q \}$, then we
  will use the notation $\ensuremath{\operatorname{Zer}} (Q,
  \R^k)$ instead. We will denote by
  $\mathcal{Q}^h$ (respectively, $Q^h$) the \emph{homogenizations} of the polynomials
  in $\mathcal{Q}$ (respectively, the polynomial $Q$), and denote by
  $\ensuremath{\operatorname{Zer}} (\mathcal{Q}^h,
  \mathbb{P}^k_{\C})$ (respectively,
  $\ensuremath{\operatorname{Zer}} (Q^h,
  \mathbb{P}^k_{\C})$) the common zeros of the
  family $\mathcal{Q}^h$ (respectively, the polynomial $Q^h$) in the
  projective space $\mathbb{P}^k_{\C}$. 
\end{notation}

\begin{notation}
  For any $Q \in \R [X_1, \ldots, X_k]$ we will
  denote by $\deg (Q)$ the degree of $Q$. More, generally for a tuple of
  polynomials $\mathcal{Q}= (Q_1, \ldots, Q_{\ell}) \in
  \R [X_1, \ldots, X_k]^{\ell}$ we will denote
  $\deg (\mathcal{Q}) = (d_1, \ldots, d_{\ell})$ where $d_i = \deg (Q_i)$, $1
  \leq i \leq \ell$.
\end{notation}

\begin{notation}
  For any semi-algebraic set $S \subset \R^k$,
  we will denote by $b_i (S)$ the $i$-th Betti number of $S$. In particular,
  $b_0 (S)$ is the number of semi-algebraically connected components of $S$. 
\end{notation}

\begin{notation}
  For any semi-algebraic set $S \subset \R^k$,
  we will denote by $\dim S$ the {{\em real dimension\/}} of $S$. For any $x
  \in S$, we denote by $\dim_x S$ the local real dimension of $S$ at $x$. Note
  that unlike complex varieties, an irreducible real variety can have have
  different local dimensions at different points.
\end{notation}

\begin{remark}
  We will at times slightly abuse notation and use the same letter to denote a
  tuple of polynomials as well as the ordered finite set whose elements are
  the elements of the tuple. This should not cause any confusion.
\end{remark}

The following theorem gives a reasonably tight bound on the number of semi-algebraically connected components of the realizations of all realizable sign conditions of a finite family of polynomials restricted to
a variety.
It generalizes earlier results of Alon {\cite{Alon}},
Warren {\cite{Warren}} and Pollack and Roy {\cite{PR}}, 
and has  found several applications in discrete geometry.

\begin{theorem}\cite{BPR8}
  \label{thm:BPR} Let $\mathcal{P}, \mathcal{Q} \subset
  \R [X_1, \ldots, X_k]$ be finite families of
  polynomials such that the degrees of the polynomials in $\mathcal{P},
  \mathcal{Q}$ are bounded by $d$,
  $\ensuremath{\operatorname{card}}\mathcal{P}= s$, and
  $\dim_{\R} (V) = k'$, where $V
  =\ensuremath{\operatorname{Zer}} (\mathcal{Q},
  \R^k)$. Then,
\[
    \sum_{\sigma \in \{ 0, 1, - 1 \}^{\mathcal{P}}} b_0
    (\ensuremath{\operatorname{Reali}} (\sigma, V))  \leq  O (1)^k s^{k'} d^k .
\]
\end{theorem}

Notice that in the bound in Theorem \ref{thm:BPR}, while the exponent of $s$
depends on the dimension of the variety $V$, the exponent of $d$ is that of
the ambient space. Moreover, the bound depends only on the maximum degree of 
the polynomials in $\mathcal{P}$ and $\mathcal{Q}$. This is a consequence of
the fact that the proof involves taking \emph{sums of squares} of the
polynomials in $\mathcal{P}$ and $\mathcal{Q}$, and thus only the maximum
degree plays a role in the argument. This feature of taking the sum of squares
is something that is common in the proofs of all the bounds mentioned above.
As such they all depend on the \emph{maximum} of the degrees of the
polynomials used to define the given set or sign conditions.

More recently, a new application of the bounds described above in discrete and computational
geometry, triggered by the work of Guth and Katz {\cite{Guth-Katz}}, raised
the question whether even the part of the bound  in Theorem \ref{thm:BPR}
that depends only on the degree $d$ could have a finer dependence on the
degrees of the polynomials in $\mathcal{P}$ and $\mathcal{Q}$, in the case
when the degrees of the polynomials in $\mathcal{Q}$ and those in
$\mathcal{P}$ differ significantly (see
{\cite{Guth-Katz, Solymosi-Tao, Matousek11, Matousek11b, zahl2011improved,
Matousek-Safernova2014}).
This is one of the primary motivations behind the results proved in the
current paper (see Section \ref{subsec:incidence} below for more detail). A second motivation is to prove a version of the Bezout
inequality on bounding the number of isolated complex solutions (or more
generally the number of connected components) of an affine polynomial system
by the product of the degrees, over real closed fields where the original
statement of the inequality does not hold (see Section \ref{subsec:bezout}, and in particular 
Example \ref{ex:fulton} and Remark \ref{rem:real-bezout} below).

A first step was taken in this direction in {\cite{Barone-Basu11a}} where the
authors of the current paper proved the following theorem (actually a more
precise statement appears in {\cite{Barone-Basu11a}} but the following
simplified version is what is important for the present purpose).

\begin{theorem}\cite{Barone-Basu11a}
  \label{thm:refined1} Let $\mathcal{P}, \mathcal{Q} \subset
  \R [X_1, \ldots, X_k]$ be finite subsets of
  polynomials such that $\deg (Q) \leq d_1$ for all $Q \in \mathcal{Q}$, $\deg
  (P) \leq d_2$ for all $P \in \mathcal{P}$. Suppose also that $d_1 \leq
  d_2$, and the real dimension of $V =\ensuremath{\operatorname{Zer}}
  (\mathcal{Q}, \R^k)$ is $k_1 \leq k$, and that
  $\ensuremath{\operatorname{card}}\mathcal{P}= s$. Then,
  \begin{eqnarray*}
    \sum_{\sigma \in \{ 0, 1, - 1 \}^{\mathcal{P}}} b_0
    (\ensuremath{\operatorname{Reali}} (\sigma, V)) & \leq & O (1)^k (s
    d_2)^{k_1} d_1^{k - k_1} .
  \end{eqnarray*}
\end{theorem}

\begin{remark}
  One should compare Theorem \ref{thm:refined1} with Theorem \ref{thm:BPR}.
  The new aspect of Theorem \ref{thm:refined1} is the more refined dependence
  on the two different degrees, taking into account the dimension of the
  variety $V$ and the fact that 
 $d_1\leq d_2$. 
\end{remark}

Notice that Theorem \ref{thm:refined1} implies the following corollary about
the number of semi-algebraically connected components of real varieties.

\begin{theorem-corollary}
  \label{cor:refined1-algebraic} Let $Q_1, Q_2 \in
  \R [X_1, \ldots, X_k]$ such that $\deg (Q_1)
  \leq d_1$,$\deg (Q_2) \leq d_2$, $d_1 \leq d_{2 }$. Let $V_1
  =\ensuremath{\operatorname{Zer}} (Q_1, \R^k)$
  and $\dim V_1 \leq k_1$, and let $V_2 =\ensuremath{\operatorname{Zer}} (\{
  Q_{1,} Q_2 \}, \R^k)$. Then,
  \begin{eqnarray*}
    b_0 (V_2) & \leq & O (1)^k d_1^{k - k_1} d_{2}^{k_1} .
  \end{eqnarray*}
\end{theorem-corollary}

\begin{proof}
 In Theorem \ref{thm:refined1}, take
$\mathcal{Q}= \{ Q_1 \}$, $\mathcal{P}= \{ Q_2 \}$ and $V =
V_1$.
\end{proof}

\subsection{Applications in discrete geometry}
\label{subsec:incidence}
While Theorem \ref{thm:refined1} (in particular, also Corollary
\ref{cor:refined1-algebraic}) has already proved useful in certain
applications in discrete and computational geometry (see
\cite{AMS12, Solymosi-Tao,Matousek-Safernova2014}), some even more recent developments
seem to require a more detailed analysis. 

The requirement of refined bounds from real algebraic geometry in the
applications mentioned above originates
in the so called ``polynomial partitioning'' method due to Guth and Katz
\cite{Guth-Katz}, which provides a framework for proving bounds in several types of problems in discrete geometry
involving finite sets of points (such as incidence problems \cite{Solymosi-Tao}, unit and distinct distance problems 
\cite{Guth-Katz, Matousek11b, zahl2011improved}
etc.). 

The original polynomial partitioning result states that given any set, $S$, of $n$ points in
$\R^k$, and an auxiliary parameter $r, \; 0 < r < n$, there exists a polynomial $P \in \R [X_1,
\ldots, X_k]$ of degree at most $O \left( r^{\frac{1}{k}} \right)$, having the
property that each semi-algebraically connected component $C$ of
$\R^k \setminus \ensuremath{\operatorname{Zer}}
(P, \R^k)$ contains at most $\frac{n}{r}$ of the
points of $S$. The number of such semi-algebraically connected components $C$
(using for instance Theorem \ref{thm:BPR}) is bounded by $O (r)$, and it is at
this point that a quantitative bound on the number of semi-algebraically
connected components of a semi-algebraic set or sign conditions enters the
proof. The polynomial partitioning theorem is a tool to decompose a given
problem involving the set $S$ into sub-problems of smaller size (corresponding
to the point sets $C \cap S$ where $C$ is a semi-algebraically connected
component of $\R^k \setminus
\ensuremath{\operatorname{Zer}} (P, \R^k)$). 
However, it might happen that most or even all the points of $S$ are
contained in $\ensuremath{\operatorname{Zer}} (P,\R^k)$
which is problematic for  a ``divide-and-conquer'' type argument. 
In this case, an obvious idea is to try
to extend the polynomial partitioning theorem to varieties of lower
dimensions, and continue the partitioning recursively. 
However, in order to prove the
strongest result possible using this approach, one 
requires tight bounds on the number of 
semi-algebraically connected components of real varieties defined by a
sequence polynomials of strictly increasing degrees,
which has a much more refined dependence on the sequence of degrees
than what was provided in Theorem \ref{thm:refined1} mentioned above
(where the length of the sequence is restricted to at most $2$).
Note however that in the  applications related to the polynomial partitioning method,  
the length of the sequence of degrees
could be as large as the dimension of the ambient space,  and Theorem \ref{thm:refined1}
is insufficient to deal with cases with
degree sequences of lengths greater than $2$.
The main result of this paper (Theorem \ref{thm:refined-algebraic}) is geared towards handling this situation, and
has already proved useful in applications involving the polynomial partitioning technique.
For example, Theorem \ref{thm:refined-algebraic} plays a crucial role in a recent application of  
multi-level polynomial partition technique for proving the tightest known bound for the point-hypersurface incidence
problem in $\mathbb{R}^4$ \cite[Theorem~1.5]{Basu-Sombra2014}.

\subsection{Failure of the naive version of Bezout inequality over the reals}
\label{subsec:bezout}
Before stating our results let us
consider what kind of refined bounds are plausible. In the case of a real
variety $V$ of $\R^k$, which is a non-singular
complete intersection (even at infinity) and defined by polynomials of degrees
$d_1 \leq d_2 \leq \cdots \leq d_{\ell}$, the number of semi-algebraically
connected components of $V$ is bounded by (see Proposition
\ref{prop:hirzebruch} as well as Remark \ref{rem:Katz} below)
\[ O (1)^k d_1 \ldots d_{\ell - 1} d_{\ell} d_{\ell}^{k - \ell} . \]
Notice that $k - \ell = \dim V$. It is thus natural to hope that such a bound
continues to hold even if the given variety is not a non-singular complete
intersection -- namely, one might hope that the number of semi-algebraically
connected components of a real variety $V \subset
\R^k$ defined by a sequence of $\ell$
polynomials having degrees $d_1 \leq d_2 \leq \cdots \leq d_{\ell}$ is bounded
by $O (1)^k d_1 \ldots d_{\ell} d_{\ell}^{\dim V}$. However, the following
well known (counter-)example (that appears in {\cite{Fulton}}) already shows
that this is not the case.

\begin{example}
  \label{ex:fulton} 
  Let $k = 3 $ and let
  \begin{eqnarray*}
    Q_1 & = & X_{3,}\\
    Q_2 & = & X_{3,}\\
    Q_3 & = & \sum_{i = 1}^2 \left( \prod^d_{j = 1} (X_i - j)^2 \right)
    .  \end{eqnarray*}

  The real variety defined by $\{ Q_1, Q_2, Q_3 \}$ is $0$-dimensional,
  and has $d^2$ isolated (in $\R^3$) points, whereas the degree sequence $(d_1, d_2,
  d_3) = (1, 1, 2 d)$, and thus the conjectured bound is $d_1 \ldots d_{\ell}
  d_{\ell}^{\dim V} = O (d)$. In particular, this example shows that the (naive version of)
  Bezout inequality which states that the number of isolated complex zeros of
  a system of polynomial equations is bounded by the product of the degrees of
  the polynomials appearing in the system, is not true over if we replace the
  complex numbers by a real closed field. Notice however that the polynomials 
  $Q_1,Q_2,Q_3$ do not define a non-singular complete intersection. 
\end{example}

While this might seem discouraging at first glance, one way to repair the
situation is to formulate a bound that depends not just on the degree sequence
and the dimension of the last variety $V = V_3
=\ensuremath{\operatorname{Zer}} (\{ Q_1, Q_2, Q_3 \},
\R^k) $, but also takes into account the
dimensions of the intermediate varieties $V_1 =\ensuremath{\operatorname{Zer}}
(Q_1, \R^k)$, $V_2
=\ensuremath{\operatorname{Zer}} (\{ Q_1, Q_2 \},
\R^k)$ etc. Notice that in Example
\ref{ex:fulton} the dimensions $k_1 = \dim V_1$, and $k_2 = \dim V_2$ are both
equal to $2$, whereas $k_3 = \dim V_3 = 0$. The number of semi-algebraically
connected components in this case is bounded by $O (d_1^{k - k_1} d_2^{k_1 -
k_2} d_3^{k_2})$, where $d_i = \deg Q_i$. This is the starting point of the
formulation of the new results proved in this paper.

We prove the following theorems where the shapes of the bounds should be seen
in the light of Example \ref{ex:fulton}.

\subsection{Main results}

For the rest of the paper we will use the following notation.
\begin{notation}
\label{not:main}
Let $Q_1, \ldots, Q_{\ell} \in \R [X_1, \ldots,
X_k]$ such that for each $i, 1 \leq i \leq \ell$, $\deg (Q_i) \leq d_i$. For
$1 \leq i \leq \ell$, denote by $\mathcal{Q}_i = \{Q_1, \ldots, Q_i \}$, $V_i
=\ensuremath{\operatorname{Zer}} (\mathcal{Q}_i,
\R^k)$, and
$\dim_{\R} (V_i) \leq k_i$. We set $V_0 =
\R^k$, and adopt the convention that $k_i = k$
for $i \leq 0$. It is clear that $k = k_0 \geq k_1 \geq \cdots \geq k_{\ell}$.
\end{notation}

With these assumptions we have the following generalization of Corollary
\ref{cor:refined1-algebraic}.

\begin{theorem}
  \label{thm:refined-algebraic}
  
Suppose that
\[ 2 \leq d_1 \leq d_2 \leq \frac{1}{k + 1} d_3 \leq \frac{1}{(k + 1)^2} d_4
   \leq \cdots \leq \frac{1}{(k + 1)^{\ell - 2}} d_{\ell} . \]
   Then,
  \begin{eqnarray*}
    b_0 (V_{\ell}) & \leq & O (1)^k \sum_{\tau = (\tau_0, \tau_1, \ldots,
    \tau_{\ell - 1})} F (k, \tau) \left( d_{\ell}^{\tau_{\ell - 1}} \prod_{1
    \leq i < \ell} ((k - \tau_{i - 1} + 1) d_i)^{\tau_{i - 1} - \tau_i}
    \right)
  \end{eqnarray*}
  where the sum on the right hand side is taken over all $\tau \in
  \mathbb{N}^{\ell}$, with $k = \tau_0 \geq \tau_1 \geq \cdots \geq
  \tau_{_{\ell - 1}} \geq 0$, and $\tau_i \leq k_i$, for each $i, 1 \leq i <
  \ell$, and
  \[ F (k, \tau) = (k - \tau_{\ell - 1} + 1) \binom{k - \tau_{\ell -
     1}}{\tau_0 - \tau_1, \tau_1 - \tau_2, \ldots, \tau_{\ell - 2} -
     \tau_{\ell - 1}} . \]
  This implies that
  \begin{eqnarray*}
    b_0 (V_{\ell}) & \leq & O (1)^{\ell} O (k)^{2 k} \left( \prod_{1 \leq j <
    \ell} d_j^{k_{j - 1} - k_j} \right) d_{\ell}^{k_{\ell - 1}} ,
  \end{eqnarray*}
  and in particular if $\ell \leq k,$
  \[ b_0 (V_{\ell}) \leq O (k)^{2 k} \left( \prod_{1 \leq j < \ell} d_j^{k_{j
     - 1} - k_j} \right) d_{\ell}^{k_{\ell - 1}} . \]
\end{theorem}

\begin{remark}
  \label{rem:complex} Note that since the real dimension of each variety $V_i$
  is at most the complex dimension of $V_i$, Theorem
  \ref{thm:refined-algebraic} remains true if we replace real dimension by
  complex dimension in the statement. This observation is important in the application of Theorem \ref{thm:refined-algebraic} to incidence problems (see \cite{Basu-Sombra2014}).
\end{remark}

\begin{remark}
  \label{rem:real-bezout} In view of Example \ref{ex:fulton} above, Theorem
  \ref{thm:refined-algebraic} can be viewed as a weak version of the Bezout
  inequality over real closed fields.
\end{remark}

The following slight modification of Example \ref{ex:fulton} shows that the
dependence on the degrees in the bound in Theorem \ref{thm:refined-algebraic}
cannot be improved.

\begin{example}
  \label{ex:tight} Let $k = k_0 \geq k_1 \geq k_2 \geq \cdots \geq k_{\ell}
  = 0$, and $d_1, \ldots, d_{\ell}$ be even. 
  For 
  $1 \leq i \leq \ell$,
  let
  \begin{eqnarray*}
    Q_i & = & \sum_{j = k - k_{_{i - 1}} + 1}^{k - k_i} \left( \prod_{h =
    1}^{d_i / 2} (X_j - h) \right)^2.
  \end{eqnarray*}
  
  Then, for $0 \leq i \leq \ell$, $\deg (Q_i) = d_i$, the real dimension
  of the variety $V_i =\ensuremath{\operatorname{Zer}} (\mathcal{Q}_i,
  \R^k)$ where $\mathcal{Q}_i = \{ Q_1, \ldots,
  Q_i \}$, is clearly $k_i$, and
  \begin{eqnarray*}
    b_0 (V_{\ell}) & = & \frac{1}{2^k} d_1^{k_0 - k_1} d_2^{k_1 - k_2} \cdots
    d_{\ell - 1}^{k_{\ell - 2} - k_{\ell - 1}} d_{\ell}^{k_{\ell - 1}} .
  \end{eqnarray*}
\end{example}

With the same assumptions as in Theorem \ref{thm:refined-algebraic}, suppose
additionally that $\mathcal{P} \subset \R [X_1,
\ldots, X_k] $ is a finite family of polynomials with $\deg (P) \leq d$ for
all $P \in \mathcal{P}$, and $\ensuremath{\operatorname{card}}\mathcal{P}=
s$, and suppose that $d_{\ell} \leq \frac{1}{k + 1} d$.

\begin{theorem}
  \label{thm:refined-semi-algebraic}
  \begin{eqnarray}
    \sum_{\sigma \in \{ 0, 1, - 1 \}^{\mathcal{P}}} b_0
    (\ensuremath{\operatorname{Reali}} (\sigma, V_{\ell})) & \leq & \sum_{j =
    0}^{k_{\ell}} 4^j \binom{s}{j} O (1)^k \Delta  \label{eqn:precise-bound}
  \end{eqnarray}
  where $\Delta$ is defined by
  \[ \Delta = \sum_{\tau = (\tau_0, \tau_1, \ldots, \tau_{\ell})} F (k, \tau)
     d^{\tau_{\ell}} \left( \prod_{1 \leq i \leq \ell} ((k - \tau_{i - 1} +
     1) d_i)^{\tau_{i - 1} - \tau_i} \right), \]
  where the sum is taken over all $\tau \in \mathbb{N}^{\ell + 1}$, with $k =
  \tau_0 \geq \tau_1 \geq \cdots \geq \tau_{_{\ell}} \geq 0$, and $\tau_i \leq
  k_i$, for each $i, 1 \leq i \leq \ell$, and
  \[ F (k, \tau) = (k - \tau_{\ell} + 1) \binom{k - \tau_{\ell}}{\tau_0 -
     \tau_1, \tau_1 - \tau_2, \ldots, \tau_{\ell - 1} - \tau_{\ell}} . \]
  This implies that
  \begin{eqnarray*}
    \sum_{\sigma \in \{ 0, 1, - 1 \}^{\mathcal{P}}} b_0
    (\ensuremath{\operatorname{Reali}} (\sigma, V_{\ell})) & \leq & O
    (1)^{\ell} O (k)^{2 k} (s d)^{k_{\ell}}  \left( \prod_{1 \leq j \leq \ell}
    d_j^{k_{j - 1} - k_j} \right) .
  \end{eqnarray*}
  In particular, if $\ell \leq k$,
  \begin{eqnarray*}
    \sum_{\sigma \in \{ 0, 1, - 1 \}^{\mathcal{P}}} b_0
    (\ensuremath{\operatorname{Reali}} (\sigma, V_{\ell})) & \leq & O (k)^{2
    k} (s d)^{k_{\ell}}  \left( \prod_{1 \leq j \leq \ell} d_j^{k_{j - 1} -
    k_j} \right) .
  \end{eqnarray*}
\end{theorem}

\begin{remark}
  Notice that in the case $\ell = 1$, the bound (\ref{eqn:precise-bound}) in
  Theorem \ref{thm:refined-semi-algebraic} implies that of Theorem
  \ref{thm:refined1}, and hence Theorem \ref{thm:refined-semi-algebraic} is a
  strict generalization of Theorem \ref{thm:refined1}.
\end{remark}

With the same assumptions as in Theorem \ref{thm:refined-semi-algebraic}, let
for $P \in \mathcal{P}$, $d_P = \deg (P)$, and for any subset $\mathcal{I}
\subset \mathcal{P}$ let
\begin{eqnarray}
\label{eqn:dsubI}
  d_{\mathcal{I}} & = & (k +
  1)^{\binom{\ensuremath{\operatorname{card}}\mathcal{I}}{2} + (k_{\ell}
  -\ensuremath{\operatorname{card}}\mathcal{I}) 
  (\ensuremath{\operatorname{card}}\mathcal{I}- 1)} \left( \prod_{P \in
  \mathcal{I}} d_P \right) (\max_{P \in \mathcal{I}} d_P)^{k_{\ell}
  -\ensuremath{\operatorname{card}}\mathcal{I}} .
\end{eqnarray}
We have the following variant of Theorem \ref{thm:refined-semi-algebraic} (the
extra precision with respect to the degrees of polynomials in $\mathcal{P}$
might be useful in applications in incidence geometry).

Using Notation \ref{not:main} and notation introduced in \eqref{eqn:dsubI} above: 

\begin{theorem}
  \label{thm:refined-semi-algebraic-2}{$\mbox{}$}
  \[
    \sum_{
    \sigma \in \{ 0, 1, - 1 \}^{\mathcal{P}}} b_0
    (\ensuremath{\operatorname{Reali}} (\sigma, V_{\ell})) 
    \leq 
    \sum_{\substack{
    \mathcal{I} \subset \mathcal{P} \\
     j =\ensuremath{\operatorname{card}}\mathcal{I} \leq k_{\ell}}}
      4^j O
    (1)^{\ell} O (k)^{2 k} d_{\mathcal{I}}  \left( \prod_{1 \leq j \leq \ell}
    d_j^{k_{j - 1} - k_j} \right).
  \]
\end{theorem}

\begin{remark}
  The condition on the degrees in Theorems \ref{thm:refined-algebraic} and
  \ref{thm:refined-semi-algebraic} might look unnatural at first glance but is
  forced on us by the method of the proof, which involves taking minors of
  matrices of size at most $(k + 1) \times (k + 1)$ with entries which are
  polynomials of degree $d_i$, $1 \leq i \leq \ell $. We want at each step,
  the degree $d_i$ to majorize the degree of the polynomial obtained as a
  minor in the previous step whose entries have degree at most $d_j$, where $j
  < i$. Notice that in the case $\ell = 2$, the condition on the degree
  sequence is just $d_1 \leq d_2$, and this allows us to recover Theorem
  \ref{thm:refined1} from Theorem \ref{thm:refined-semi-algebraic}.
\end{remark}

\begin{remark}
  We also note that in {\cite{RV96}} the authors define the
  ``complexification'' of a semi-algebraic set as the smallest complex variety
  containing it, and prove an effective bound on the geometric degree of this
  complexification which depend amongst other quantities on the real dimension
  of the given set. This degree could be thought of as the ``real degree'' of
  the semi-algebraic set. It is possible that Theorem
  \ref{thm:refined-semi-algebraic} could serve as an alternative basis for a
  good definition of the ``real degree'' of a real variety -- in the sense
  that the ``real degree'' of a real variety $V$ should control the number of
  semi-algebraically connected components of the intersection of $V$ with any
  real hypersurface of sufficiently large degree. We do not pursue this idea
  further in the current paper.
\end{remark}

Finally, we conjecture that the bounds in Theorems \ref{thm:refined-algebraic}
and \ref{thm:refined-semi-algebraic} extend to the sum of all the Betti
numbers (instead of just the zero-th one). The techniques developed in this
paper are not sufficient to prove this conjecture.

\subsection{Outline of the 
proofs of the main theorems}
\label{subsec:outline}
The main difficulty that one faces in order to prove bounds having the shapes
of Theorems \ref{thm:refined-algebraic} and Theorem
\ref{thm:refined-semi-algebraic} is that in order to respect the degree
sequence one has to be careful about taking ``sums of squares'' which spoil
the dependence on the degrees. The crucial idea is to use the notion of
``approximating'' varieties. An approximating variety is a variety which is
infinitesimally close to the given variety of the same dimension, but having
good algebraic properties which allow one to give a precise bound on the
number of its semi-algebraically connected components in terms of the sequence
of degrees of polynomials defining it (rather than just the maximum degree).
If the given variety can be covered (in a technical sense made precise later)
by a small number of such approximating varieties, then the problem of
bounding the number of semi-algebraically connected components of the given
varieties reduces to the problem of bounding the total number of
semi-algebraically connected components of these approximating varieties.

The idea of using approximating varieties originates in algorithmic
semi-algebraic geometry and it was used in {\cite{BPR95b}} to give efficient
algorithms for computing sample points on varieties and in {\cite{BPR99}} to
compute roadmaps of semi-algebraic sets. The combinatorial part of the
complexities of these algorithms depends on the dimension of the given variety
rather than that of the ambient space, and this is where the approximating
varieties play an important role in those papers. In quantitative
semi-algebraic geometry, the notion of approximating varieties was used in
{\cite{Barone-Basu11a}} in order to prove Theorem \ref{thm:refined1}.

The approximation scheme that we use, which is a generalization of the one
used in {\cite{Barone-Basu11a}} is described in Section
\ref{subsec:construction-of-approximating-varieties}  below. One difficulty in generalizing the scheme
in {\cite{Barone-Basu11a}} is that the non-singularity of polar varieties of
smooth hypersurfaces with respect to generic projections that is used in that
paper no longer holds for smooth varieties of higher co-dimension. A second
difficulty is that the sequence of local (real) dimensions at a point $x \in
V_{\ell}$ of the varieties $V_1, \ldots, V_{\ell}$ is not globally constant,
but is only a local invariant (see Example \ref{eg:main} and Figure \ref{fig:example4d}). Thus, one cannot expect to have a single global
approximating variety with good properties. We overcome the latter problem by
taking into account all possible sequences of local dimensions whether they
actually occur or not (indexed by the set $A$ below), and construct
approximating varieties with acceptable degree sequences to approximate each
of them. 

Consider the subset of points of $U_i$ of $V_{\ell}$ having local dimension
$i \leq k_{\ell}$. At each point $x \in U_i$ the dimension of $V_{\ell - 1}$
is between $i$ and $k_{_{\ell - 1}}$. Suppose we have already constructed
approximations of subsets of $V_{\ell - 1}$ consisting points having some
fixed local dimension at $V_{\ell - 1}$. Using these approximations and adding
appropriately many equations in each case we construct a set of approximations
of $U_i$. Taking all these approximating varieties, for all $i, 0 \leq i \leq
k_{\ell}$, and noticing that $V_{\ell}$ is the union of the $U_i$'s we obtain
a global approximation of $V_{\ell}$ (see Example \ref{eg:prop:main} and Figures \ref{fig:V}, \ref{fig:V2}, \ref{fig:V1},
\ref{fig:V11} and \ref{fig:V21} below).

More precisely, we construct a family of basic semi-algebraic sets each of
the form,
\begin{eqnarray*}
  \Bas(\mathcal{P}, \mathcal{Q}) & : = & \{ x \in
  \R'^k \mid P (x) = 0 , Q (x) \leq 0, P \in
  \mathcal{P}, Q \in \mathcal{Q} \},
\end{eqnarray*}
where $\R'$ is some real closed extension of
$\R$ depending on the particular approximating
set. The family of pairs $\{ (\mathcal{P}^{\alpha}_{\tau, \ell},
\mathcal{Q}^{\alpha}_{\tau, \ell}) \}_{\tau \in A \subset \mathbb{N}^{\ell},
\alpha \in I (\tau)}$ defining these approximating varieties are indexed by a
pair of indices $\tau, \alpha$ coming from two finite set of indices $A_\ell 
\subset \mathbb{N}^{\ell}$, and $I_{\ell} (\tau)$. While the definition of the
second, $I_{\ell} (\tau)$, is a bit technical and which we defer for later,
the definition of the index set $A_\ell$ is the following.
\begin{eqnarray*}
  A_\ell & := & \{ \tau = (\tau_1, \ldots, \tau_{\ell}) \in \mathbb{N}^{\ell} \mid
  k \geq \tau_1 \geq \tau_2 \geq \cdots \geq \tau_{\ell} , \tau_i \leq k_i, 1
  \leq i < \ell \} .
\end{eqnarray*}
For any given $\tau \in A_\ell$, let $V_{\tau} \subset V_{\ell}$ denote the closure
of the set of points $x \in V_{\ell}$ such that the local real dimension of
$V_i$ at $x$ is equal to $\tau_i$, for each $i, 1 \leq i \leq \ell$. The union
of the approximating sets $V^{\alpha}_{\sigma, \ell} =
\Bas(\mathcal{P}^{\alpha}_{\sigma}, \mathcal{Q}^{\alpha}_{\sigma})$ with $\sigma
\leq \tau$, ``approximates'' $V_{\tau}$ in a certain precise sense (see
Proposition \ref{prop:main} below), and since clearly $V_{\ell} =
\bigcup_{\tau \in A} V_{\tau}$, the union of all the approximating sets $\{
V^{\alpha}_{\tau, \ell} \}_{\tau \in A \subset \mathbb{N}^{\ell}, \alpha \in
I (\tau)}$ approximate the whole variety $V_{\ell}$. Because of the
approximating property, in order to bound the number of semi-algebraically
connected components of $V_{\ell}$ it suffices to bound the sum of the number
of semi-connected components of each one of the approximating sets
$V^{\alpha}_{\tau, \ell}$. The tuples $\mathcal{P}^{\alpha}_{\tau, \ell},
\mathcal{Q}^{\alpha}_{\tau, \ell}$ have the following properties that enable
us to obtain good bounds on the number of semi-algebraically connected
components of $V^{\alpha}_{\tau, \ell}$ (see Proposition
\ref{prop:transversal} below).
\begin{enumerate}[a{\textup{)}}]
  \item The tuple of polynomials $\mathcal{P}^{\alpha}_{\tau, \ell}$ define a
  non-singular, bounded complete intersection of dimension $\tau_{\ell} \leq
  k_{\ell}$. In particular, this means that the cardinality of
  $\mathcal{P}^{\alpha}_{\tau, \ell}$ is equal to $k - \tau_{\ell}$. Suppose
  that $\mathcal{P}^{\alpha}_{\tau, \ell} = (P_1, \ldots, P_{k -
  \tau_{\ell}})$. Let for $1 \leq i \leq \ell$, $\ell_i = \tau_{i - 1} -
  \tau_i$, with the convention that $\tau_0 = k$, and $L_i = \sum_{h = 1}^i
  \ell_h$. Then for each $i , 1 \leq i \leq \ell$, the degrees of the
  polynomials $P_{L_{i - 1} + 1}, \ldots, P_{L_i}$ are bounded by $O (k d_i)$.
  
  \item $\mathcal{Q}^{\alpha}_{\tau, \ell}$ is either empty or contains one
  polynomial, $Q^{\alpha}_{\tau, \ell}$, with $\deg (Q^{\alpha}_{\tau, \ell})
  = O (d_{\ell})$, and $\mathcal{P}', Q^{\alpha}_{\tau, \ell}$, where
  $\mathcal{P}'$ is any subset of $\mathcal{P}^{\alpha}_{\tau, \ell}$, defines
  a non-singular complete intersection.
\end{enumerate}
It remains to bound the number of semi-algebraically connected components of
each $V^{\alpha}_{\tau, \ell}$ and take the sum of these bounds, for which
we use the same result as in {\cite{Barone-Basu11a}} where a bound is derived
using a classical formula for the Betti numbers of complex non-singular
complete intersections and the Smith inequality (see Proposition
\ref{prop:hirzebruch} below). The number of approximating varieties (which is
independent of the given degree sequence) and the bounds on the degree
sequences of their defining polynomials as stated in Properties a) and b)
above are good enough to give us the bound in Theorem
\ref{thm:refined-algebraic}.

Theorem \ref{thm:refined-semi-algebraic} follows from Theorem
\ref{thm:refined-algebraic} using standard techniques already used in
{\cite{BPR95a}} and no fundamentally new ingredients.

The rest of the paper is organized as follows. In Section
\ref{sec:preliminaries}, we recall some basic facts about real closed fields
of Puiseux series that we need for making deformation arguments. We also
recall some results proved in {\cite{Barone-Basu11a}} on the choice of generic
coordinates. Finally, in Section \ref{sec:proofs} we prove the main theorems.

\section{Preliminary results}\label{sec:preliminaries}

\subsection{Deformation of several equations to general
position}\label{subsec:generalposition}
In this section we describe how to deform a system of equations using infinitesimals so that
the set of common zeros of the deformed equations (in certain real closed non-archimedean extensions
of the ground field)  has good properties. For this
we first need to recall some properties of Puiseux series with coefficients in a real
closed field. We refer the reader to {\cite{BPRbook2}} for further detail.

We begin with some notation.

\begin{notation}
  For $\R$ a real closed field we denote by
  $\R \langle \varepsilon \rangle$ the real
  closed field of algebraic Puiseux series in $\varepsilon$ with coefficients
  in $\R$. We use the notation
  $\R \langle \varepsilon_1, \ldots,
  \varepsilon_m \rangle$ to denote the real closed field
  $\R \langle \varepsilon_1 \rangle \langle
  \varepsilon_2 \rangle \cdots \langle \varepsilon_m \rangle$. Note that in
  the unique ordering of the field $\R \langle
  \varepsilon_1, \ldots, \varepsilon_m \rangle$, $0 < \varepsilon_m \ll
  \varepsilon_{m - 1} \ll \cdots \ll \varepsilon_1 \ll 1$. Also, note that
  both fields $\R \langle \varepsilon \rangle,
  \R \langle \delta \rangle$ are sub-fields in a
  natural way of $\R \langle \varepsilon, \delta
  \rangle$. 
\end{notation}

\begin{notation}
\label{not:ext}
  If $\R'$ is a real closed extension of a real
  closed field $\R$, and $S \subset
  \R^k$ is a semi-algebraic set defined by a
  first-order formula with coefficients in $\R$,
  then we will denote by $\ensuremath{\operatorname{Ext}} (S,
  \R') \subset
  \R'^k$ the semi-algebraic subset of
  $\R'^k$ defined by the same formula. It is
  well-known that $\ensuremath{\operatorname{Ext}} (S,
  \R')$ does not depend on the choice of the
  formula defining $S$ {\cite{BPRbook2}}.
\end{notation}

\begin{notation}
  For $x \in \R^k$ and $r \in
  \R$, $r > 0$, we will denote by $B_k (x, r)$
  the open Euclidean ball centered at $x$ of radius $r$. If
  $\R'$ is a real closed extension of the real
  closed field $\R$ and when the context is
  clear, we will continue to denote by $B_k (x, r)$ the extension
  $\ensuremath{\operatorname{Ext}} (B_k (x, r),
  \R')$. This should not cause any confusion.
\end{notation}

\begin{notation}
  For elements $x \in \R \langle \varepsilon
  \rangle$ which are bounded over $\R$ we denote
  by $\lim_{\varepsilon} x$ to be the image in
  $\R$ under the usual map that sets
  $\varepsilon$ to $0$ in the Puiseux series $x$.
\end{notation}

\begin{notation}
  \label{not:deformation} Let $Q \in \R [X_1,
  \ldots, X_k]$, $0 \leq q \leq k$, and $H \in
  \R [X_{q + 1}, \ldots, X_k]$. Let $\zeta$ be a
  new variable. We denote
  \begin{eqnarray*}
    \ensuremath{\operatorname{Def}} (Q, \zeta, q, H) & = & (1 - \zeta) Q -
    \zeta H.  \label{eqn:not:Def}
  \end{eqnarray*}
\end{notation}

\begin{notation}
  \label{not:def} For $\mathcal{P} = (P_1, \ldots, P_m) $, with each $P_i \in
  \R [X_1, \ldots, X_k]$, $1 \leq q \leq k$, and
  $\mathcal{G}= (G_1, \ldots, G_m)$ with each $G_i \in
  \R [X_{q + 1}, \ldots, X_k]$, and $\zeta$ a
  new variable, we denote by $\ensuremath{\operatorname{Def}} (\mathcal{P},
  \zeta, q, \mathcal{G})$ the tuple
  \[ (\ensuremath{\operatorname{Def}} (P_1, \zeta, q, G_1), \ldots,
     \ensuremath{\operatorname{Def}} (P_m, \zeta, q, G_m)) , \]
  and by $\ensuremath{\operatorname{Def}} (\mathcal{P}, \zeta, q,
  \mathcal{G})^h$ the corresponding tuple of homogenized polynomials
  \[ (\ensuremath{\operatorname{Def}} (P_1, \zeta, q, G_1)^h, \ldots,
     \ensuremath{\operatorname{Def}} (P_m, \zeta, q, G_m)^h)  . \]
\end{notation}

\begin{notation}
  \label{not:jacobian} For $\mathcal{F}= (F_1, \ldots, F_{k - p}) , q \leq p
  \leq k $, we denote the jacobian matrix
  \begin{eqnarray*}
    \mathcal{\ensuremath{\operatorname{Jac}}} (\mathcal{F}, p, q) & : = &
    \left(\begin{array}{ccc}
      \frac{\partial F_1}{\partial X_{q + 1}} &  \cdots & \frac{\partial
      F_{k - p}}{\partial X_{q + 1}}\\
      \vdots &  & \vdots\\
      \frac{\partial_{F 1}}{\partial X_k} & \cdots & \frac{\partial F_{k -
      p}}{\partial X_k}
    \end{array}\right)
  \end{eqnarray*}
  whose rows are indexed by $[q + 1, k]$ and columns by $[1, k - p]$.
  
  For $J \subset [q + 1, k]$, $\ensuremath{\operatorname{card}}J = k - p$ and
  $k \in J$, let $\mathcal{\ensuremath{\operatorname{Jac}}}_J$ denote the
  $(k - p) \times (k - p)$ matrix extracted from the matrix
  $\mathcal{\ensuremath{\operatorname{Jac}}} (\mathcal{\mathcal{F}}, p, q)$ by
  extracting the rows whose index are in $J$, and let
  \[ \ensuremath{\operatorname{jac}}_J = \det
     (\mathcal{\ensuremath{\operatorname{Jac}}}_J) . \]
  Let
  \begin{eqnarray}\label{eqn:FJ}
    \mathcal{\mathcal{F}}_J & := & \mathcal{\mathcal{F}} \cup \bigcup_{i \in
    [q + 1, k] \setminus J} \{ \ensuremath{\operatorname{jac}}_{J \cup \{ i \}
    \setminus \{ k \}} \},
  \end{eqnarray}
  and the finite constructible set
  \begin{eqnarray}\label{eqn:CJ}
    C_J (\mathcal{F}) & := & \{ x \in \ensuremath{\operatorname{Zer}} (
    \mathcal{\mathcal{F}}_J, \R^k)  \mid
    \ensuremath{\operatorname{jac}}_J (x) \neq 0 \} .
  \end{eqnarray}
\end{notation}

\begin{proposition}
  \label{prop:simple} Let $\mathcal{F}= (F_1, \ldots, F_{k - p})$, each $F_i
  \in \R [X_1, \ldots, X_k]$, and such that the
  variety $\ensuremath{\operatorname{Zer}} (\mathcal{F}^h,
  \mathbb{P}^k_{\C})$ is a non-singular complete
  intersection. Let $x \in \R^k$ be a
  non-generate critical point of the projection map to the $X_k$-coordinate
  restricted to the variety $V =\ensuremath{\operatorname{Zer}} (\mathcal{F},
  \R^k)$. Then, there exists a subset $J \subset
  [1, k]$, $\ensuremath{\operatorname{card}}J = k - p$, $k \in J$, satisfying
  the following two conditions.
  \begin{enumerate}[1.]
    \item The $(k - p) \times (k - p)$ matrix, $\tmop{Jac}_J$, extracted from
    the matrix $\mathcal{\tmop{Jac}} (\mathcal{\mathcal{F}}, p, 0)$ by
    extracting the rows whose index are in $J$, evaluated at $x$ is
    non-singular.
    
    \item The point $x$ is a simple zero of the system $\mathcal{\mathcal{F}}_J$ (see \eqref{eqn:FJ} for definition).
   \end{enumerate}
\end{proposition}

\begin{proof}
First note that using the Jacobian criteria for non-singularity of real algebraic varieties
(see for example \cite[Definition 3.3.4]{BCR}), we have that
the variety $V$ is of dimension equal to $p$ and non-singular.
Moreover, $x$ is a critical point of the projection map to the $X_k$
coordinate restricted to $V$, by the inverse function theorem we can choose
$p$ coordinates (not including $X_k$) such that the remaining $k - p$
co-ordinates of points of $V$ in a small enough neighborhood $U$ of $x$ are
smooth functions of these chosen $p$ co-ordinates. Without loss of generality
let these $p$ coordinates be $X_1, \ldots, X_p$. We will denote the remaining
co-ordinate functions on $U$ by $X_{p + 1} (X_1, \ldots, X_p), \ldots, X_k
(X_1, \ldots, X_p)$ noting that they are smooth semi-algebraic functions of
$X_1, \ldots, X_p$.

We use that
\begin{enumerate}
  \item $\tmop{Jac} (\mathcal{F}, p, 0) (x)$ has full rank since $x$ is a
  non-singular point of $V$, and
  
  \item $\tmop{Hess} (X_k (X_1, \ldots, X_p)) (x)$ is non-singular since $x$
  is a non-degenerate critical point with respect to $X_k$. 
\end{enumerate}
Let $J = [p + 1, k]$, and consider the Jacobian matrix $\tmop{Jac}
(\mathcal{F}_J, 0, 0)$.
\[ \tmop{Jac} (\mathcal{F}_J, 0, 0) = \left(\begin{array}{cccccc}
     \frac{\partial F_1}{\partial X_1} & \ldots & \frac{\partial F_{k -
     p}}{\partial X_1} &  \frac{\partial \tmop{jac}_{J \cup \{ 1 \} \setminus
     \{ k \}}}{\partial X_1} & \ldots  & \frac{\partial
     \tmop{jac}_{J \cup \{ p \} \setminus \{ k \}}}{\partial X_1}\\
     \vdots &  & \vdots & \vdots &  & \vdots\\
     \frac{\partial F_1}{\partial X_k} & \ldots &  \frac{\partial F_{k -
     p}}{\partial X_k} &  \frac{\partial \tmop{jac}_{J \cup \{ 1 \} \setminus
     \{ k \}}}{\partial X_k} & \ldots & \frac{\partial \tmop{jac}_{J \cup \{ p
     \} \setminus \{ k \}}}{\partial X_k}
   \end{array}\right) .\]
Since, by definition of the functions $X_{p + 1} (X_1, \ldots, X_p), \ldots,
X_k (X_1, \ldots, X_p)$
\[ F_i (X_1, \ldots, X_p, X_{p + 1} (X_1, \ldots, X_p), \ldots, X_k (X_1,
   \ldots, X_p)) \equiv 0, \]
for $1 \leq i \leq k - p$, by the chain rule for $1 \leq j \leq p,$
\begin{eqnarray}
\label{eqn:F} 
  0 &=& \frac{\partial F_1}{\partial X_j} + \frac{\partial F_1}{\partial X_{p +
  1}}  \frac{\partial X_{p + 1}}{\partial X_j} + \ldots + \frac{\partial
  F_1}{\partial X_k}  \frac{\partial X_k}{\partial X_j},   \nonumber\\
  & \vdots &   \nonumber\\
   \\
  0 &=& \frac{\partial F_{k - p}}{\partial X_j} + \frac{\partial F_{k -
  p}}{\partial X_{p + 1}}  \frac{\partial X_{p + 1}}{\partial X_j} + \ldots +
  \frac{\partial F_{k - p}}{\partial X_k}  \frac{\partial X_k}{\partial X_j} .
    \nonumber
\end{eqnarray}
Let $\Delta = \det \tmop{Jac} (\mathcal{F}, p, p)$. Notice that in the
sub-matrix $\tmop{Jac} (\mathcal{F}, p, 0)$ of $\tmop{Jac} (\mathcal{F}_J, 0,
0)$, for each $1 \leq j \leq p  $, adding
\[ \sum_{i = p + 1}^k \frac{\partial X_i}{\partial X_j} \cdot \tmop{row}_i
   (\tmop{Jac} (\mathcal{F}, p, 0)) \]
to the $j$-th row and using \eqref{eqn:F} we can clear out the first $p$ rows.
Since, $\tmop{rank} (\tmop{Jac} (\mathcal{F}, p, 0) (x)) = k - p$, this
implies that $\Delta (x) \neq 0$.

From Cramer's Rule, we have
\begin{eqnarray*}
  \frac{\partial X_k}{\partial X_1}  & = & \frac{- \tmop{jac}_{J \cup \{ 1 \}
  \setminus \{ k \}}}{\Delta},\\
  & \vdots & \\
  \frac{\partial X_k}{\partial X_p} & = & \frac{- \tmop{jac}_{J \cup \{ p \}
  \setminus \{ k \}}}{\Delta} .
\end{eqnarray*}
Let for $1 \leq i \leq p$, \
\[ G_i (X_1, \ldots, X_p) = - \tmop{jac}_{J \cup \{ i \} \setminus \{ k \}}
   (X_1, \ldots, X_p, X_{p + 1} (X_1, \ldots, X_p), \ldots, X_k (X_1, \ldots,
   X_p))  . \]
Substituting above we get that

\begin{eqnarray*}
  \frac{\partial X_k}{\partial X_1} & = & \frac{G_1 (X_1, \ldots,
  X_p)}{\Delta},\\
  & \vdots & \\
  \frac{\partial X_k}{\partial X_p} & = & \frac{G_p (X_1, \ldots,
  X_p)}{\Delta} .
\end{eqnarray*}
From the quotient rule,
\begin{eqnarray*}
  \frac{\partial^2 X_k}{\partial X_i \partial X_j} & = & \frac{\frac{\partial
  G_i}{\partial X_j} \Delta - G_i  \frac{\partial \Delta}{\partial X_j}}{\Delta^2},
\end{eqnarray*}
and in particular
\[ \tmop{Hess} (X_k) (x) = \left(\begin{array}{c}
     \frac{\partial^2 X_k}{\partial X_i \partial X_j} (x)
   \end{array}\right)_{1 \leqslant i, j \leqslant p} = \left(\begin{array}{c}
     \frac{\frac{\partial G_i}{\partial X_j} (x)}{\Delta (x)}
   \end{array}\right)_{1 \leqslant i, j, \leqslant p} \]
noticing that since $x$ is a critical point of the function $X_k$ restricted
to $V$, $G_1 (x) = \cdots = G_p (x)$ = 0.

Applying the chain rule again we have that for $1 \leq i, j \leq p$,
\begin{equation}
  \frac{\partial G_i}{\partial X_j} = - \frac{\partial \tmop{jac}_{J \cup
  \{ i \} \setminus \{ k \}}}{\partial X_j} - \frac{\partial \tmop{jac}_{J
  \cup \{ i \} \setminus \{ k \}}}{\partial X_{p + 1}}  \frac{\partial X_{p +
  1}}{\partial X_j} - \cdots - \frac{\partial \tmop{jac}_{J \cup \{ i \}
  \setminus \{ k \}}}{\partial X_k}  \frac{\partial X_k}{\partial X_j}.
  \label{eqn:G}
\end{equation}
Finally, for each $1 \leq j \leq p  $, adding
\[ \sum_{i = p + 1}^k \frac{\partial X_i}{\partial X_j} \cdot \tmop{row}_i
   (\tmop{Jac} (\mathcal{F}_J, 0, 0)) \]
to the $j$-th row,  
and using \eqref{eqn:F} and \eqref{eqn:G}, we see that
$\tmop{Jac} (\mathcal{F}_J, 0, 0) (x)$ is row equivalent to the matrix
\[ \left(\begin{array}{cc}
     \tmmathbf{0} & - \frac{\tmop{Hess} (X_k)}{\Delta} (x)\\
     \mathbf{I}_{k - p} & \ast
   \end{array}\right) \]
which is clearly non-singular, since $x$ is a non-degenerate critical point of
$X_k$, which implies that the $\tmop{Hess} (X_k) (x)$ is non-singular, and we
have already observed that $\Delta (x) \neq 0$.
\end{proof}

\begin{definition}
  Let $X \subset \mathbb{P}^k_{\C}$ be a
  non-singular variety, and $(H_{\mu})_{\mu = (\mu_0 : \mu_1)\in \PP^1_\C}$ a pencil of
  hyperplanes. We call the pencil of varieties $(X_{\mu} = X \cap
  H_{\mu})_{\mu}$ a \emph{Lefschetz pencil} if it satisfies the two following
  conditions.
  \begin{enumerate}[1.]
    \item The base locus $B$ is smooth of co-dimension two in $X$.
    
    \item Each member $X_{\mu}$ of the pencil has at most one ordinary double
    point as a singularity.
  \end{enumerate}
  
\end{definition}

The main result about Lefschetz pencil we will require is the following well known result from complex algebraic geometry (see for example  \cite[Corollary 2.10]{Voisin2}).

\begin{proposition}
  \label{prop:Lefschetz} If $X \subset
  \mathbb{P}^k_{\C}$ is a non-singular variety,
  then any generic pencil of hyperplane sections of $X$ is Lefschetz.
\end{proposition}

\begin{remark}
  \label{rem:generic} Observe that a generic tuple of polynomials
  $\mathcal{G}= (G_1, \ldots, G_{k - p})$ where each $G_i \in
  \R [X_1, \ldots, X_k]$ with $\deg (H_i) = d_i$
  and is chosen generically, will have the property that the variety $W
  =\ensuremath{\operatorname{Zer}} (\mathcal{G},
  \mathbb{P}^k_{\C})$ is non-singular and the
  pencil of hyperplane sections $(W_{\mu} = W \cap H_{\mu})_{\mu}$ indexed by
  $\mu = (\mu_0 : \mu_1)$, where $H_{\mu} \subset
  \mathbb{P}_{\C}^k$ is defined by the equation
  $\mu_0 X_0 + \mu_1 X_k = 0$, is a Lefschetz pencil for the variety $W$ by
  Proposition \ref{prop:Lefschetz} above. 
\end{remark}

Let $0 \leq p \leq k$, 
$\bar\eps= (\eps_1,\ldots,\eps_m)$ be a tuple of variables,
and $\mathcal{P} = (P_1, \ldots, P_{k - p})$, $P_i \in
\R\la\bar\eps\ra [X_1, \ldots, X_k]$ with $\deg P_i \leq d_i$,
and $P \in \R [X_1, \ldots, X_k]$, $\deg P \leq
d$. Let $0 \leq q < p \leq k$ and $\mathcal{G}= (G_1, \ldots, G_{k - p})$ be
a tuple of polynomials with $G_i \in \R [X_{q
+ 1}, \ldots, X_k]$ with $\deg (G_i) = d_i$, and $G \in
\R [X_{q + 1}, \ldots, X_k]$ be another
polynomial with $\deg (G) = d$, such that
\begin{enumerate}[1.]
  \item The variety $W = \tmop{Zer} (\mathcal{G} \cup \{G\}, \mathbb{P}^{k -
  q}_\C)$ is a non-singular complete
  intersection.
  
  \item The pencil of hyperplane sections $(W_{\mu} = W \cap H_{\mu})_{\mu}$
  indexed by $\mu = (\mu_0 : \mu_1)\in \PP^1_\C$, where $H_{\mu} \subset
  \mathbb{P}_\C^{k - q}$ is defined by the
  equation $\mu_0 X_0 + \mu_1 X_k = 0$, is a Lefschetz pencil for the variety
  $W$. 
\end{enumerate}

\hide{
Let for every $z \in \R^q$,
\begin{eqnarray}\label{eqn:Fz}
\mathcal{F}_z &=& \ensuremath{\operatorname{Def}} (\mathcal{P}, \zeta, q,
   \mathcal{G}) (z, \cdot), \ensuremath{\operatorname{Def}} (P, \delta,
   q, G) (z, \cdot) .
   \end{eqnarray}
   }
   
We also need the following notation.

\begin{notation}
  \label{not:projection} For $1 \leq p \leq q \leq k$, we denote by $\pi_{[p,
  q]}$ the projection map on the coordinates $X_p, \ldots, X_q$, and also
  denote by $\R^{[p, q]}$ the subspace spanned
  by these coordinates. For any set $S \subset
  \R^k$, and $z \in
  \R^{[1, p]}$ we will denote by $S_z$ the fiber
  $S \cap \pi^{- 1}_{[1, p]} (z)$.
\end{notation}

\begin{proposition}
  \label{prop:transversal} 
 For every $z \in
  \R^q$, the following holds.
  \begin{enumerate}[1.]
    \item \label{item:prop:transversal1}
    $\tmop{Def} (\mathcal{P}, \eta, q, \mathcal{G}) (z, \cdot)^h,
    \tmop{Def} (P, \delta, q, G) (z, \cdot)^h$ defines a non-singular
    complete intersection $V_z \subset
    \mathbb{P}_{\C \langle \delta, \bar\eps,\eta
    \rangle}^{k - q}$ of dimension $p - q - 1$.
    
    \item \label{item:prop:transversal2} The pencil of hyperplane sections $(V_{z, \mu} = V_z \cap
    H_{\mu})_{\mu}$ indexed by $\mu = (\mu_0 : \mu_1)$, where $H_{\mu} \subset
    \mathbb{P}_{\C \langle \delta, \bar\eps,\eta
    \rangle}^{k - q}$ is defined by the equation $\mu_0 X_0 + \mu_1 X_k = 0$,
    is a Lefschetz pencil for the variety $V_z$.
    
    \item \label{item:prop:transversal3}For each singular point $x \in \C^k$ of
    the pencil $(V_{z, \mu})_{\mu}$, there exists $J \subset [k - q + 1, k]$,
    $\tmop{card} J = k - p$ and $k \in J$, such that $x \in C_J(\mathcal{F}_z)$ (see \eqref{eqn:CJ} 
    for definition),
     and $x$ is a simple zero of the system $(\mathcal{F}_{z})_J$ (see \eqref{eqn:FJ} for definition),
      where   
\begin{eqnarray*}
\mathcal{F}_z &=& \ensuremath{\operatorname{Def}} (\mathcal{P}, \eta, q,
   \mathcal{G}) (z, \cdot), \ensuremath{\operatorname{Def}} (P, \delta,
   q, G) (z, \cdot) .
   \end{eqnarray*}
    
  \end{enumerate}
\end{proposition}

\begin{proof}
Replacing $\eta$ and $\delta$ by new variables
$s$ and $t$ (respectively), and setting $s = t = 1$ we have that $(\tmop{Def}
(\mathcal{P}, 1, q, \mathcal{G}) (z, \cdot)^h, \tmop{Def} (P, 1, q, G) (z,
\cdot)^h) = (\mathcal{G}^h, G^h)$ define a non-singular complete intersection
in $W \subset \mathbb{P}_{\C\la\bar\eps\ra}^{k - q}$ (by
hypothesis). Moreover, the pencil of hyperplane sections $(W_{\mu} = W \cap
H_{\mu})_{\mu}$ is Lefschetz by hypothesis. Since the property of being a
non-singular complete intersection as well as a fixed pencil of hyperplane
section being Lefschetz is stable, it also holds for an open neighborhood of
the point $(s, t) = (1, 1)$. The set of pairs $(s, t)$ for which any of these
two properties is violated is Zariski closed, defined over $\C\la\bar\eps\ra$, and is not the whole of
$\mathbb{P}^1_{\C\la\bar\eps\ra} \times
\mathbb{P}^1_{\C\la\bar\eps\ra}$. In particular the
complement contains 
the point $(\eta,\delta)$ (since $\eta,\delta$ are algebraically independent  over $\C\la\bar\eps\ra$).
  This proves parts
\ref{item:prop:transversal1}. and \ref{item:prop:transversal2}. of the proposition. 
Part \ref{item:prop:transversal3}.  follows from Proposition
\ref{prop:simple}.
\end{proof}

We also need the following proposition.

\begin{proposition}
  \label{prop:algebraic} Let $C$ be a bounded s.a. connected component of
  $\Bas(\mathcal{P}, \mathcal{Q})$. Then, there exists a subset subset
  $\mathcal{Q}' \subset \mathcal{Q}$, and a semi-algebraically connected
  component $D$ of $\ensuremath{\operatorname{Zer}} (\mathcal{P} \cup
  \mathcal{Q}', \R^k)$ such that $D \subset C$.
\end{proposition}

\begin{proof}
 See Proposition 13.1 in
{\cite{BPRbook2}}.
\end{proof}

\begin{proposition}
  \label{prop:simple2} Let $\mathcal{F}= (F_1, \ldots, F_k)$ be a tuple of
  polynomials with $F_i \in \R [X_1, \ldots,
  X_k]$ and let $\mathcal{G}= (G_1, \ldots, G_k)$, with $G_i \in
  \R [X_1, \ldots, X_k]$, be a 
   tuple of
  polynomials with $\deg (\mathcal{G}) \leq \deg (\mathcal{F})$. Let $x \in
  \R^k$ a simple zero of $\mathcal{F}$. Then,
  there exists a simple zero $\tilde{x} \in \ensuremath{\operatorname{Zer}}
  (\ensuremath{\operatorname{Def}} (\mathcal{F}, \zeta, 0, \mathcal{G}),
  \R \langle \zeta \rangle^k)$, such that
  $\lim_{\zeta} \tilde{x} = x$.
\end{proposition}

\begin{proof}
 It follows from the fact that $x$ is a simple zero
of the family $\mathcal{F}$ that any infinitesimal perturbation of the family
$\mathcal{F}$ will have a simple zero, $\tilde{x} \in
\C \langle \zeta \rangle^k$, in an infinitesimal
neighborhood of $x$. 
To see this observe that since $x \in \R^k$ (and is thus in particular bounded over $\R$), 
it belongs to the image under the $\lim_\zeta$ map (extended to elements of $\C\la \zeta\ra^k$ which are bounded over $\R$) of
$\ensuremath{\operatorname{Zer}}
  (\ensuremath{\operatorname{Def}} (\mathcal{F}, \zeta, 0, \mathcal{G}),
  \C \langle \zeta \rangle^k)$. Thus, there exists
 $\tilde{x} \in  \ensuremath{\operatorname{Zer}}
  (\ensuremath{\operatorname{Def}} (\mathcal{F}, \zeta, 0, \mathcal{G}),
  \C \langle \zeta \rangle^k)$, such that
  $\lim_{\zeta} \tilde{x} = x$.
 Moreover, since $x$ is  a simple zero of $\ZZ(\mathcal{F},\R^k)$, we have that 
 $\det(\mathrm{Jac}(\mathcal{F},0,0))(x) \neq 0$ (see Notation \ref{not:jacobian}). 
 Since, 
 \[\lim_\zeta(\det(\mathrm{Jac}
  (\ensuremath{\operatorname{Def}} (\mathcal{F}, \zeta, 0, \mathcal{G}),0,0))(\tilde{x})) = 
  \det(\mathrm{Jac}(\mathcal{F},0,0))(x),
  \]
   this implies that 
   \[
   \det(\mathrm{Jac}(
  (\ensuremath{\operatorname{Def}} (\mathcal{F}, \zeta, 0, \mathcal{G}),0,0)))(\tilde{x})
  \neq 0
  \]
  as well, 
  since $\det(\mathrm{Jac}(\mathcal{F},0,0))(x) \in \R \setminus \{0\}$,
  and hence
   $\tilde{x}$ is a simple zero of $\ensuremath{\operatorname{Def}} (\mathcal{F}, \zeta, 0, \mathcal{G})$.
   
Moreover, $\tilde{x}$ must belong to
$\R \langle \zeta \rangle^k$ as long as the
perturbed polynomials also have real coefficients. Otherwise, since complex
zeros must occur in conjugate pairs, if $\tilde{x} \not\in
\R \langle \zeta \rangle^k$, then $\tilde{x} \neq
\overline{\tilde{x}}$, while $\lim_{\zeta} \tilde{x} = \lim_{\zeta}
\overline{\tilde{x}} = x$, and this implies that $x$ is not a simple zero of
$\mathcal{F}$.
\end{proof}

\subsection{Generic coordinates}

We recall in this section a result proved in {\cite{Barone-Basu11a}} that we
will require.

\begin{notation}
  For a real algebraic set $V =\ensuremath{\operatorname{Zer}} (Q,
  \R^k) $we let $\ensuremath{\operatorname{reg}}
  (V)$ denote the non-singular points in dimension $\dim V$ of $V$ (Definition
  3.3.9 in {\cite{BCR}}).
\end{notation}

\begin{definition}
  Let $V =\ensuremath{\operatorname{Zer}} (Q,
  \R^k) $ be a real algebraic set. Define $V^k =
  V$, and for $0 \leq i \leq k - 1$ define
  \[ V^{(i)} = V^{(i + 1)} \setminus \ensuremath{\operatorname{reg}} (V^{(i +
     1)}) . \]
  Let $d_V (i)$ denote the dimension of $V^{(i)}$.
\end{definition}

\begin{definition}
  Let $V =\ensuremath{\operatorname{Zer}} (Q,
  \R^k)$ be a real algebraic set, $1 \leq j \leq
  k$, and $\ell \in \ensuremath{\operatorname{Gr}} (k, k - j)$. We say that
  the linear space $\ell$ is $j$-good with respect to $V$ if either:
  \begin{itemize}
    \renewcommand{\labelitemi}{$\bullet$}\renewcommand{\labelitemii}{$\bullet$}\renewcommand{\labelitemiii}{$\bullet$}\renewcommand{\labelitemiv}{$\bullet$}\item
    $j \not\in d_V ([0, k])$,
    
    \item or $d_V (i) = j$, and the set
    \[ 
    \{ x \in \tmop{reg} (V^{(i)}) \mid \tmop{dimT}_x
       V^{(i)} \cap \ell = 0 \}  \]
    is a non-empty dense Zariski open subset of $V^{(i)}$.
  \end{itemize}
\end{definition}

\begin{definition}
  Let $V =\ensuremath{\operatorname{Zer}} (Q,
  \R^k)$ and $B = \{ v_1, \ldots, v_k \}$ be a
  basis of $\R^k$. We say that the basis $B$
  is \emph{generic} with respect to $V$ if for each $j, 1 \leq j \leq k$,
  the linear space $\ensuremath{\operatorname{span}} (v_1, \ldots, v_{k - j})$
  is $j$-good with respect to $V$. {\color{red} {\textsc{{\color{red}
  {\color{red} }}}}}
\end{definition}

The following proposition appears in {\cite{Barone-Basu11a}}.

\begin{proposition}
  \label{prop:generic} Let $V =\ensuremath{\operatorname{Zer}} (Q,
  \R^k)$ and $\{ v_1, \ldots, v_k \}$ be a
  basis of $\R^k$. Then, there exists a
  non-empty open semi-algebraic subset $U$ of linear transformations
  $\ensuremath{\operatorname{GL}} (k, \R)$ such
  that for every $T \in U$ the basis $\{ T (v_1), \ldots, T (v_k) \}$ is
  generic with respect to $V$. 
\end{proposition}

\section{Proofs of the main theorems}\label{sec:proofs}

We now fix polynomials $Q_1, \ldots, Q_{\ell}$ and and the varieties $V_1,
\ldots, V_{\ell}$ as in Theorem \ref{thm:refined-algebraic}. We will assume
if necessary by initially squaring each polynomial that each $Q_i$ is
non-negative over $\R^k$. Since this increases
each degree by a multiplicative factor of $2$, this does not affect the
asymptotics of the bound.

The section is organized as follows. In Subsection \ref{subsec:construction-of-approximating-varieties} we define
certain approximating semi-algebraic sets and prove their important properties. In Subsection
\ref{subsec:bounds-on-Betti} we recall and then apply in the current context  certain well-known 
bounds on the Betti numbers of non-singular cmplete intersections. Finally, we prove the main theorems of this
paper in Subsections \ref{subsec:proof-of-theorem-algebraic} and \ref{subsec:proof-of-theorem-semi-algebraic}.

\subsection{Definitions and main properties of approximating semi-algebraic sets}
\label{subsec:construction-of-approximating-varieties}
We first introduce in \ref{subsubsec:notation} some necessary notation, and then in Subsection
\ref{subsec:approx} below we describe the construction of certain
semi-algebraic sets approximating the varieties $V_j$. The main properties of
these sets is then proved in Subsection \ref{subsec:properties}. The
approximating properties of these sets are  proved in Proposition
\ref{prop:main},  and the quantitative estimates on the degrees of the
polynomials appearing in the description of these approximating sets is proved
in Proposition \ref{prop:quantitative}.

\subsubsection{Notation}
\label{subsubsec:notation}
\begin{notation}
  \label{not:tauofx}
   For any semi-algebraic set $S$ and $x \in S$, we denote
  by $\dim_x S$ the local dimension of $S$ at $x$. For $0 \leq j \leq \ell$
  and $x \in V_j$, we denote
  \[ \dim^{(j)} (x) = (\dim_x V_1, \ldots, \dim_x V_j) . \]
\end{notation}

\begin{notation}
\label{not:poset}
We will use the natural partial order on the sets $\mathbb{N}^j$, and denote for
$\sigma = (\sigma_1,\ldots,\sigma_j),\tau=(\tau_1,\ldots,\tau_j) \in \mathbb{N}^j$, 
$\sigma \leq \tau$, if 
$\sigma_i \leq \tau_i$ for all $1 \leq i \leq j$.
\end{notation}

Before proceeding further we illustrate the notation introduced above by considering some examples. 

We first consider again Example \ref{ex:fulton} from Section \ref{sec:intro}.

In this example, (following Notation \ref{not:main})  we have
\begin{eqnarray*}
V_1 &=& \ZZ(X_3,\R^k), \\
V_2 &=&  \ZZ(X_3,\R^k),\\
V_3 &=&  \{1,\ldots,d\}^3.
\end{eqnarray*}

The various functions $\dim^{(j)}:V_j \rightarrow \mathbb{N}^j, j=1,2,3$ (cf. Notation
\ref{not:tauofx}) are as follows.

\begin{eqnarray*}
\dim^{(1)}(x) &=& (2) \mbox{ if } x\in V_1, \\
 \dim^{(2)}(x) &=& (2,2) \mbox{ if } x\in V_2=V_1,\\
\dim^{(3)}(x) &=& (2,0) \mbox{ if } x\in V_3.
\end{eqnarray*}

The next example is slightly more involved but is helpful in understanding the
proof of Proposition \ref{prop:main} below.

\begin{example}
\label{eg:main}
Let $k=4$, $\ell = 3$, and 
\begin{eqnarray*}
Q_1  &=& (X_1+X_2+X_3+X_4)(X_3^2+X_4^2),\\
Q_2  &=&  X_3^2 + X_4^2, \\
Q_3  &=& X_3^2 + X_4^2 + (X_1^2- X_2^3)^2.
\end{eqnarray*}

We denote by $(e_i)_{1\leq i \leq 4}$ the elementary basis vectors in $\R^4$.
Denote by $L_1$ the hyperplane defined by $X_1+X_2+X_3+X_4=0$, by
$L_2=\mathrm{span}(e_1,e_2)$, 
the linear subspace defined by $X_3=X_4=0$, and by $C_3$ the cubic curve
contained in $L_2$, defined by the equation $X_1^2 - X_2^3=0$. 
Notice that $L_1 \cap L_2$ is the line in $\mathrm{span}(e_1,e_2)$ defined by $X_1+X_2=0$, and it meets $C_3$ at the points $x^0 = \mathbf{0}$ (which is a singular point of $C_3$), and
$x^1 = (-1,1,0,0,0)$ (which is a regular point of $C_3$).  
These sets are depicted in Figure \ref{fig:example4d}.

\begin{figure}
\includegraphics[scale=0.5]{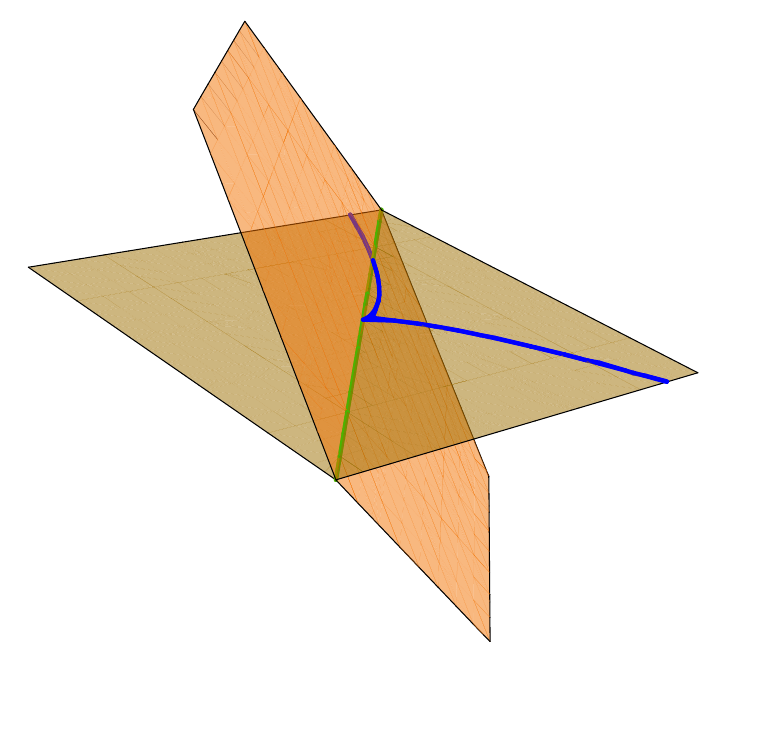}
\caption{The sets $L_1,L_2$ and $C_3$ in Example \ref{eg:main} restricted to the hyperplane $X_4=0$}
\label{fig:example4d}
\end{figure}

We have using Notation \ref{not:main},
\begin{eqnarray*}
V_1 &=& L_1 \cup L_2, \\
V_2 &=& L_2, \\
V_3 &=& C_3.
\end{eqnarray*}

The various functions $\dim^{(j)}:V_j \rightarrow \mathbb{N}^j, j=1,2,3$ can now be described as follows.

\begin{eqnarray*}
\dim^{(1)}(x) &=& (3) \mbox{ if } x\in L_1, \\
 \dim^{(1)}(x) &=& (2) \mbox{ if } x\in L_2\setminus L_1,\\
\dim^{(2)}(x) &=& (3,2) \mbox{ if } x\in L_2 \cap L_1, \\
 \dim^{(2)}(x) &=& (2,2) \mbox{ if } x\in L_2 \setminus L_1,\\
\dim^{(3)}(x) &=& (3,2,1) \mbox{ if } x = x^0,x^1, \\
 \dim^{(3)}(x) &=& (2,2,1) \mbox{ if } x\in C_3\setminus \{x^0,x^1\}. 
\end{eqnarray*}
\end{example}

\begin{remark}
\label{rem:ex-main}
Observe that in Example \ref{eg:main} above, 
the point $x^0 = \mathbf{0}$ is a singular point of $V_3$, and 
$x^0 \in V_{(3,2,1)}$. However, for any open neighborhood $U\subset V_3$ of $x_0$ in $V_3$,
and $x' \in U \setminus\{x^0\}$, we have that $x'$ is a regular point of $V_3$, and 
if moreover $x' \neq x^1$, $\dim^{(3)}(x') = (2,2,1) < (3,2,1)$ (cf. Notation \ref{not:poset}).
Notice also that $(2,2,1)$ is a minimal element of the set $\{\dim^{(3)}(x) \mid x \in V_3\}$, and
the semi-algebraic subset of $V_3$ defined by 
\[
\{x \in V_3 \mid \dim^{(3)}(x) = (2,2,1) \} = 
V_3 \setminus \{x^0, x^1\}
\]
is open in $V_3$ (cf. Proposition \ref{prop:uppersc} below).
\end{remark}

The following property of the function $\dim^{(j)}:V_j \rightarrow \mathbb{N}^j$ will be important
later.

Following Notation \ref{not:main} as before we have the following proposition.
\begin{proposition}
\label{prop:uppersc}
Let $\sigma \in \mathbb{N}^j, 1 \leq j \leq \ell$. Then, the semi-algebraic subset 
$V^{\leq \sigma}_j \subset V_j$ defined by 
\[
V^{\leq \sigma}_j = \{ x \in V_j \;\mid\; \dim^{(j)}(x) \leq \sigma \}
\]
is open in $V_j$. In particular, if  
$U \subset V_j$ is an open semi-algebraic subset of $V_j$, and 
$\sigma \in \mathbb{N}^j$ is such that 
$\sigma$  is a minimal element of the set 
$\{ \dim^{(j)}(x) \mid x \in U\}$,
then the semi-algebraic subset $
\{x \in U \mid \dim^{(j)}(x) = \sigma \}$
is open in $V_j$. 
\end{proposition}

\begin{proof}
The proof is by induction on $j$. If $j=1$, then the proposition follows immediately from the
upper semi-continuity property of the dimension function. Now suppose that the proposition is
true for all smaller values of $j$. Let $\sigma = (\sigma_1,\ldots,\sigma_j)$ and let 
$\sigma' = (\sigma_1,\ldots,\sigma_{j-1})$. Using the induction hypothesis we have that 
$V_{j-1}^{\leq \sigma'}$ is open in $V_{j-1}$. This implies that there exists an open semi-algebraic subset  $U \subset \R^k$, such that $V_{j-1}^{\leq \sigma'} = V_{j-1} \cap U$. 
Also, the  
semi-algebraic set $V^{\leq \sigma_j}_j = \{x \in V_j\;\mid\; \dim_x(V_j) \leq \sigma_j\} \subset V_j$ 
is open in $V_j$.
Thus, there exists an open semi-algebraic set $U' \subset \R^k$, such that 
$V^{\leq \sigma_j} = V_j \cap U'$. 
Now, 
\begin{eqnarray*}
V_j^{\leq \sigma} &=& V^{\leq \sigma_j}_j \cap  V_{j-1}^{\leq \sigma'} \\
                                  &=& (V_j \cap U') \cap (V_{j-1} \cap U) \\
                                  &=& V_j \cap U \cap U' \quad (\mbox{since } V_j \subset V_{j-1}).
\end{eqnarray*}
Hence, $V^{\leq \sigma}_j $ is open in $V_j$.
\end{proof}

\begin{notation}
  \label{not:tau} For $0 \leq j \leq \ell$ we call $\tau = (\tau_1, \ldots,
  \tau_j) \in \mathbb{N}^j$ {{\em admissible}} if it satisfies the
  following two conditions.
  \begin{enumerate}[1.]
    \item $\tau_1 \geq \cdots \geq \tau_j$,
    
    \item for $1 \leq i < j$, $\tau_i \leq k_i$.
  \end{enumerate}
  We denote the subset of admissible tuples of $\mathbb{N}^j$ by $A_j$, and
  denote by $A$ the set $A_{\ell}$. For $\sigma = (\sigma_1, \ldots,
  \sigma_j), \tau = (\tau_1, \ldots, \tau_j) \in A_j$, we say $\sigma \leq
  \tau$, if $\sigma_i \leq \tau_i$ for each $i, 1 \leq i \leq j$.
  \end{notation}

\begin{notation}
  \label{not:multiplepuiseux} For each $j, 1 \leq j \leq \ell$, we denote by
  $\R_j$ the real closed field
  \[
  \R \langle \delta_j, \ldots, \delta_1,
 \eta_1, \zeta_1, \ldots, \eta_j, \zeta_j \rangle.
  \]
  Notice that
  $\R_j$ is a real closed extension of the field
  $\R_{j - 1}$. For any semi-algebraic subset $S
  \subset \R_j$, we will denote by $S_b$ the
  union of semi-algebraically connected components of $S$ which are bounded
  over $\R .$
\end{notation}

\begin{remark}
  \label{rem:ordering-of-infinitesimals} For readers familiar with arguments
  in real algebraic geometry involving multiple infinitesimals, this ordering
  of the infinitesimals in Notation \ref{not:multiplepuiseux} might seem
  somewhat counter-intuitive, since we will consider the varieties $V_i$'s in
  the order $V_1$, $V_2$, etc., and the infinitesimal $\delta_i$ will be used
  to perturb the variety $V_i$, one would expect that the infinitesimals
  $\delta_i$'s to be ordered the other way round. The reason behind this
  ordering of the infinitesimals will become clear in the proof of Proposition
  \ref{prop:main} below.
\end{remark}

\subsubsection{Definition of sequences of approximating semi-algebraic
sets}
\label{subsec:approx}

We now describe the construction of certain semi-algebraic sets approximating
the varieties $V_j$. 
We  assume that $V_1$, and hence each $V_j$, are
bounded over $\R$.

We will also use the following notation. 
\begin{notation}
\label{not:binom}
For any set $X$ and $j \geq 0$ we will denote by $\binom{X}{j}$ the set of all subsets of $X$ of cardinality $j$.
\end{notation}

\begin{definition}
\label{def:tuple}  
For any $\tau \in A_j$ we define an index set $I_j (\tau)$, and a family
$(V_{\tau, j}^{\alpha} \subset \R_j^k)_{\alpha
\in I_j (\tau)}$ as follows. Each $V^{\alpha}_{\tau, j} = \left(
\Bas(\mathcal{P}^{\alpha}_{\tau, j}, \{{Q}^{\alpha}_{\tau, j}\})
\right)_b$, where $\mathcal{P}^{\alpha}_{\tau, j}$,
is an ordered tuple of polynomials,
and ${Q}^{\alpha}_{\tau, j} \in \R_j[X_1,\ldots,X_k]$
defined inductively as follows.
\begin{enumerate}[1.]
  \item If $j = 0$, then for $\tau = ()$, define $I_0 (\tau) = \{ - 1 \}$, and
  $\mathcal{P}^{(- 1)}_{\tau, 0} = (0) , {Q}^{(- 1)}_{\tau, 0}
  = 0$.
  
  \item Otherwise, we denote by $\tau' = (\tau_1, \ldots, \tau_{j - 1})$ and
  let $p = \tau_{j - 1}$, $q = \tau_j $. Let $G$ be a generic
  polynomial in $\R [X_{q + 1}, \ldots, X_k]$
  strictly positive over $\R^{k - q}$ with 
  $\deg(G) = \deg(\bar{P}_j)$,
  \begin{eqnarray*}
    \bar{P}_j  & = & 
    \sum_{1 \leq i \leq j} Q_i \in
    \R [X_{1,} \ldots, X_k],\\
    \tilde{P}_j  & = & \tmop{Def} (\bar{P}_j, \delta_j, q, H) \in
    \R \langle \delta_j \rangle [X_{1},\ldots, X_k].
  \end{eqnarray*}
  \item
  \begin{eqnarray*}
  I_j (\tau) & = & I_{j - 1} (\tau') \times \{-1\} , \: \tmop{if} \tau_{j -1} = \tau_j,\\
    & = & I_{j - 1} (\tau') \times \binom{[\tau_j + 1, k]}{k - \tau_{j - 1} +
    1}, \; \tmop{else}
  \end{eqnarray*}
  (where $\times$ denotes the usual Cartesian product).
  \item For each triple $(\alpha \in I_{j - 1} (\tau'),
  \mathcal{P}^{\alpha}_{\tau', j - 1}, {Q}^{\alpha}_{\tau', j -
  1})$
  \begin{itemize}
    \item if $\tau_{j - 1} = \tau_j$, then denoting $\beta = (\alpha, - 1)$
    let
    \begin{eqnarray*}
      \mathcal{P}^{\beta}_{\tau, j} & = & \mathcal{P}^{\alpha}_{\tau', j -
      1},\\
      {Q}^{\beta}_{\tau, j} & = & \tilde{P}_j .
    \end{eqnarray*}
    \item otherwise, suppose that
    \[ \mathcal{P}=\mathcal{P}^{\alpha}_{\tau', j - 1} = (P_1, \ldots, P_{k -
       p}) \subset \R_{j - 1}  [X_1,
       \ldots, X_k]^{k - p}, \]
    with $\deg (P_i) = d_i'$, for $1 \leq i \leq k - p$, and $\overline{d'} =
    (d_1', \ldots, d'_{k - p})$. Let
    \[ \mathcal{G}= (G_1, \ldots, G_{k - p}) \]
    be generic polynomials in $\R [X_{q + 1},
    \ldots, X_k]$ with $\deg (G_i) = d_i'$ and strictly positive over
    $\R^{k - q}$, $1 \leq i \leq k - p$.
    
    We define (using Notation \ref{not:def})
    \begin{eqnarray}
      \tilde{\mathcal{P}} & = & \tmop{Def} (\mathcal{P}, \eta_j, q,
      \mathcal{G}) \nonumber\\
      \mathcal{F} & = & (\tilde{\mathcal{P}}, \tilde{P}_j) .  \label{eqn:FF}
    \end{eqnarray}
    Finally, for each $J \in \binom{[\tau_j + 1, k]}{k - \tau_{j - 1} + 1}$,
    denoting $\beta = (\alpha, J)$, and following the notation introduced
    above (and using Notation \ref{not:jacobian})
    \begin{eqnarray}
      \mathcal{P}^{\beta}_{\tau, j} & = & \tmop{Def} (\mathcal{F}_J, \zeta_j,
      k, \mathcal{G}'),  \label{eqn:deformedF}\\
      {Q}^{\beta}_{\tau, j} & = & {Q}^{\alpha}_{\tau', j -
      1} , \nonumber
    \end{eqnarray}
    where $\mathcal{G}' = (G_1', \ldots, G'_{k - q})$ is another tuple of
    generic polynomials strictly positive over
    $\R^k$ with $\deg (\mathcal{G}') =
    (\bar{d}_{\alpha}, d_j, d', \ldots, d')$, where $d' = (k - p + 1) d_j$ and
    $\bar{d}_{\alpha} = \deg (\mathcal{P})$.
  \end{itemize}
\end{enumerate}
\end{definition}

\begin{notation}
\label{not:vtau}
For each $j, 0 \leq j \leq \ell $, $\tau \in A_j$, let (cf. Notation
\ref{not:tauofx})
\begin{eqnarray*}
  \tilde{V}_{\tau} & = & \bigcup_{\alpha \in I_j (\tau)} V_{\tau,
  j}^{\alpha},\\
  V_{\tau} & = & \overline{\{ x \in V_j \mid \dim^{(j)} (x) = \tau \}} .
\end{eqnarray*}
\end{notation}

Using Notation \ref{not:vtau}:
\begin{proposition}
  \label{prop:cover} For each $j$, $1 \leq j \leq \ell$,
  \[ V_j = \bigcup_{\tau = (\tau_1, \ldots, \tau_j) \in A_j, \tau_j \leq k_j}
     V_{\tau} . \]
\end{proposition}

\begin{proof}
This is immediate from the definition of $A_j$ and
the various $V_{\tau}$,  
\[
\tau = (\tau_1, \ldots, \tau_j) \in A_j, \tau_j \leq
k_j,
\]
and the fact that $\dim V_i \leq k_i$ for $0 \leq i \leq
j$.
\end{proof}

\subsubsection{Properties of the approximating sets}\label{subsec:properties}

The following proposition and its corollary guarantees the approximating
properties of the sets $V^{\alpha}_{\tau, j}$ defined above and is the main
technical proposition of the paper.

Assume that the given system of coordinates is generic with respect to the
finite number of varieties $V_{\tau}$ (cf. Proposition \ref{prop:generic}).

\begin{proposition}
  \label{prop:main} For all $\tau = (\tau_1, \ldots, \tau_j) \in A_j$, with
  $\tau_j \leq k_j$,
  \[ V_{\tau} \subset W_{\tau} \subset V_j , \]
  where
  \[ 
  W_{\tau} = \bigcup_{\sigma} \lim_{\delta_j} \tilde{V}_{\sigma}, 
  \]
  and the union is taken over all $\sigma \in A_j$ with $\sigma_j = \tau_j$,
  and $\sigma_i \leq \tau_i$ for all $1 \leq i < j$. 
\end{proposition}

In the proof of Proposition \ref{prop:main} we need the following technical
lemma that we prove first. We draw the attention of the reader to the ordering
of the infinitesimals in this lemma, which is particularly delicate and plays
a very important role in the proof of the lemma.
In particular, notice that if $\bar\eps=(\eps_1,\ldots,\eps_m),\delta$ are variables, then we
have the following diagram of real closed subfields of the real closed field $\R\la\delta,\bar\eps\ra$.

\[
\xymatrix{
&\R\la\delta,\bar\eps\ra &\\
\R\la\delta\ra \ar@{^{(}->}[ru]&& \R\la\bar\eps\ra\ar@{^{(}->}[lu]
\\
&\R\ar@{^{(}->}[lu]\ar@{^{(}->}[ru]\ar@{^{(}->}[uu]&
}
\]

In particular, if $V$ (respectively, $Z$) is a semi-algebraic subset of $\R\la\bar\eps\ra^k$ (respectively,
$\R\la\delta\ra^k$), then $\Ext(V,\R\la\delta,\bar\eps\ra),\Ext(Z,\R\la\delta,\bar\eps\ra)$ 
(recall Notation \ref{not:ext})
are semi-algebraic subsets of $\R\la\delta,\bar\eps\ra^k$.

\begin{lemma}
  \label{lem:technical} Let $P, H \in \R [X_1,
  \ldots, X_k]$, P non-negative, and $H$ strictly positive at all points of
  $\R^k$. Let $V \subset
  \R \langle \bar{\varepsilon} \rangle^k$ be a
  semi-algebraic set bounded over $\R$, where
  $\bar{\varepsilon} = (\varepsilon_1, \ldots, \varepsilon_m)$. Let $\tilde{P}
  = (1 - \delta) P - \delta H$, and $Z$ a semi-algebraically connected
  component of $\ensuremath{\operatorname{Zer}} (\tilde{P},
  \R \langle \delta \rangle^k)$, such that $Z
  =\ensuremath{\operatorname{Zer}} (\tilde{P},
  \R \langle \delta \rangle^k) \cap B_k (x, r)$,
  for some $x \in \R^k$ and $r > 0$, $r \in
  \R$ . Suppose that $\lim_{\varepsilon_1}
  (\ensuremath{\operatorname{Ext}} (V, \R
  \langle \delta, \bar{\varepsilon} \rangle)) \cap Z \neq \emptyset$. Then,
  $\ensuremath{\operatorname{Ext}} (V, \R
  \langle \delta, \bar{\varepsilon} \rangle) \cap
  \ensuremath{\operatorname{Ext}} (Z, \R \langle
  \delta, \bar{\varepsilon} \rangle) \neq \emptyset$.
\end{lemma}

Before proving  Lemma  \ref{lem:technical} we illustrate it with a simple example.
\begin{example}
\label{eg:technical}
In this example, $k=2$, and
\begin{eqnarray*}
V &=& \ZZ((X_1^2 + X_2^2 - 1)^2 -\eps,\R\la\eps\ra^2) \text{  (shown in blue in Figure \ref{fig:technical})}, \\  
P &=& (X_1 - 1)^2 + X_2^2, \\
H &=& 1.
\end{eqnarray*}
We display the various sets occurring in this example in Figure \ref{fig:technical} after choosing $\eps,\delta$,
to be certain sufficiently small positive real numbers,  with $\eps \ll \delta$.
It is clear from definition that $P$ is non-negative, 
$\ZZ(P,\R^k)$ consists of the single point $(1,0)$ (shown in black in Figure  \ref{fig:technical}), 
and the variety 
$\ZZ(\tilde{P},\R\la\delta\ra^2)$,
where $\tilde{P} = (1-\delta)P - \delta H$, 
has one semi-algebraically connected component,  $Z$, which is 
also depicted in black.
The semi-algebraic set
$\lim_\eps V = \ZZ((X_1^2+X_2^2-1)^2,\R^2)$ is the unit circle centered at the origin (shown in red),
and $\Ext(V,\R\la\delta,\eps\ra)$ meets $\Ext(Z,\R\la\delta,\eps\ra)$ in two points, and
$\Ext(V,\R\la\delta,\eps\ra) \cap \Ext(Z,\R\la\delta,\eps\ra)$ is not empty,  and
consists of four points
as can be seen in  Figure \ref{fig:technical}.

\begin{figure}
\includegraphics[scale=0.40]{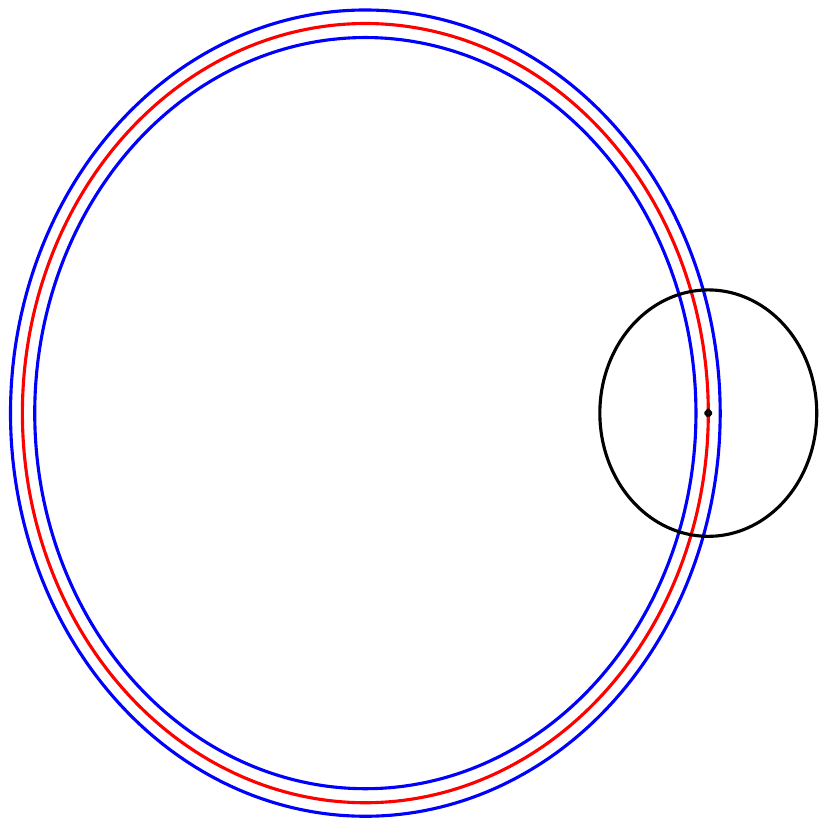}
\vspace{-5cm}
\caption{Example illustrating Lemma \ref{lem:technical}}
\label{fig:technical}
\end{figure}
\end{example}

\begin{proof}[Proof of Lemma \ref{lem:technical}]
 Let $G \in \R (X_1,
\ldots, X_k)$ denote the rational function $\frac{P}{P + H}$ which is
continuous, and takes non-negative values at all points of
$\R^k$ by hypothesis. Let $y \in \tmop{Ext} (V,
\R \langle \delta, \bar{\varepsilon} \rangle)$ be
such that $z = \lim_{\varepsilon_1} y \in Z$. Since, $Z$ is contained in
$\overline{B_k (x, r)}$, and $y \in \tmop{Ext} (V,
\R \langle \delta, \bar{\varepsilon} \rangle)$ is
$\varepsilon_1$-infinitesimally close to $z \in Z$, it is clear that
$\tmop{Ext} (V, \R \langle \delta,
\bar{\varepsilon} \rangle) \cap \overline{B_k (x, r)}$ contains $y$ and in
particular is not empty. Let $C$ be the semi-algebraically connected component
of $\tmop{Ext} (V, \R \langle \delta,
\bar{\varepsilon} \rangle) \cap \overline{B_k (x, r)}$ which contains $y$.

We prove that $\tmop{Ext} (C, \R \langle \delta,
\bar{\varepsilon} \rangle) \cap \tmop{Ext} (Z, \R
\langle \delta, \bar{\varepsilon} \rangle) \neq \emptyset$. Suppose otherwise.
Then, $G (y) \neq \delta$. Suppose without loss of generality that $G (y) -
\delta > 0$. Since, $z = \lim_{\varepsilon_1} y \in \tmop{Zer} (\tilde{P},
\R \langle \delta \rangle^k)$, it is clear that
$\lim_{\varepsilon_1} (G (y) - \delta) = 0$. Let $h = \inf_{x \in C} G (x)$.
Since, $C$ is a semi-algebraic set defined over
$\R \langle \bar{\varepsilon} \rangle$, and $G$ is
a continuous rational function defined over $\R$,
it follows that $h \in \R \langle
\bar{\varepsilon} \rangle$. Moreover, since $\tmop{Ext} (C,
\R \langle \delta, \bar{\varepsilon} \rangle) \cap
\tmop{Zer} (\tilde{P}, \R \langle \delta,
\bar{\varepsilon} \rangle^k) = \emptyset$, $G (y) - \delta > 0$, and $C$ is
closed and bounded, the infimum of $G$ over $C$ is achieved at a point, and
hence $h > \delta$. On the other hand, from the fact that 
$\lim_{\varepsilon_1} (G (y) - \delta) = 0$, it follows that
$\lim_{\varepsilon_1} h = \delta$. This is impossible, since 
$\lim_{\varepsilon_1} h \in \R$.
\end{proof}

\begin{lemma}
\label{lem:main2}
Suppose that $\sigma=(\sigma_1,\ldots,\sigma_j),\tau = (\tau_1,\ldots,\tau_j) \in A_j$ 
with $\sigma_i \leq \tau_i, 1\leq i < j$, and $\sigma_j = \tau_j$. Then, 
$W_{\sigma} \subset W_{\tau}$.
\end{lemma}
\begin{proof}
Obvious from the definitions of $W_{\sigma}$ and $W_{\tau}$.
\end{proof}

Before giving the proof of Proposition \ref{prop:main} we consider the following three-dimensional
example which illustrates some of the finer points.

\begin{example}
\label{eg:prop:main}
Let $k=3$, $\ell = 3$, and 
\begin{eqnarray*}
Q_1  &=& (X_1^2 + X_2^2 + X_3^2 -1) (X_3^2 + (X_1^2 + \frac{1}{2}X_2^2 -1)^2), \\
Q_2  &=&  (X_3^2 + (X_1^2 + \frac{1}{2}X_2^2 -1)^2),\\
Q_3  &=& X_3^2 + X_2^2 + (X_1-1)^2.
\end{eqnarray*}

The variety $V_1$ (shown in Figure \ref{fig:V}) is  bounded, and equal to the union of the unit sphere $S\subset \R^3$ (shown in orange),  and an ellipse, $C$ (shown in green), contained in the plane $\mathrm{span}(e_1,e_2)$, with $S \cap C = \{(\pm 1,0,0)\}$ (shown in red).
The variety $V_2 = C$, and $V_3 = \{(1,0,0)\}$.

\begin{figure}
\includegraphics[scale=0.75]{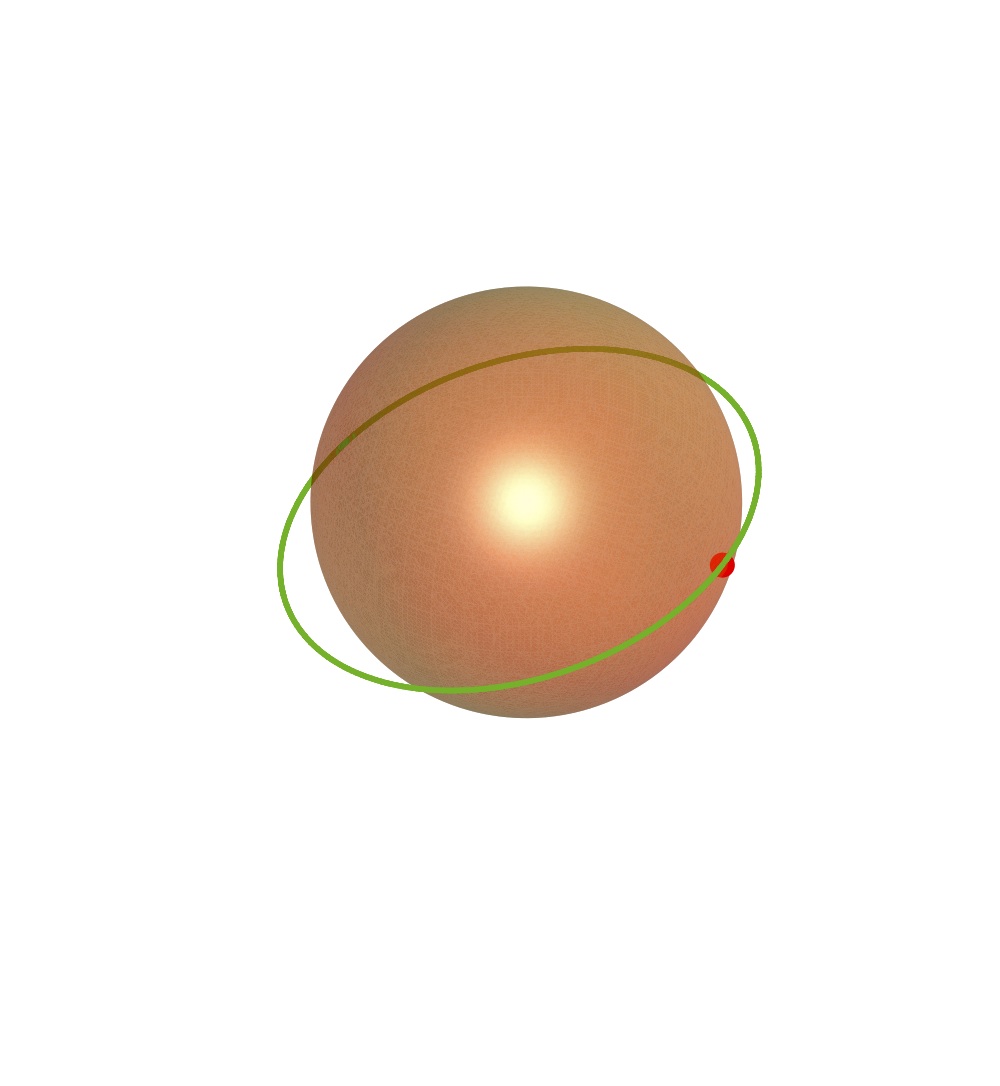}
 \vspace{-3cm}
 \caption{The variety $V_1$.}
  \label{fig:V}
 \end{figure}

The various functions $\dim^{(j)}:V_j \rightarrow \mathbb{N}^j$ are as follows.

\begin{eqnarray*}
\dim^{(1)}(x) &=& (2) \mbox{ if } x\in S, \\
 \dim^{(1)}(x) &=& (1) \mbox{ if } x\in C\setminus S,\\
\dim^{(2)}(x) &=& (2,1) \mbox{ if } x\in S\cap C, \\
 \dim^{(2)}(x) &=& (1,1) \mbox{ if } x\in C\setminus S,\\
\dim^{(3)}(x) &=& (2,1,0) \mbox{ if } x = (1,0,0).
 \end{eqnarray*}
 
 It follows that
 \begin{eqnarray*}
 V_{(2)} &=& S, \\
 V_{(1)} &=& C, \\
 V_{(1,1)} &=& C,\\
 V_{(2,1)} &=& S \cap C.
 \end{eqnarray*}

 In the Figures \ref{fig:V2}, \ref{fig:V1} and \ref{fig:V11}  below we depict the approximating semi-algebraic sets $\tilde{V}_\tau$, 
 for $\tau= (2),(1)$ and $(1,1)$, respectively. 
 Note that in order to be able to draw these pictures we used (small) finite values of the infinitesimals, and so the pictures are for illustrative purposes only.
 
 \begin{figure}[h!]
\includegraphics[scale=1.0]{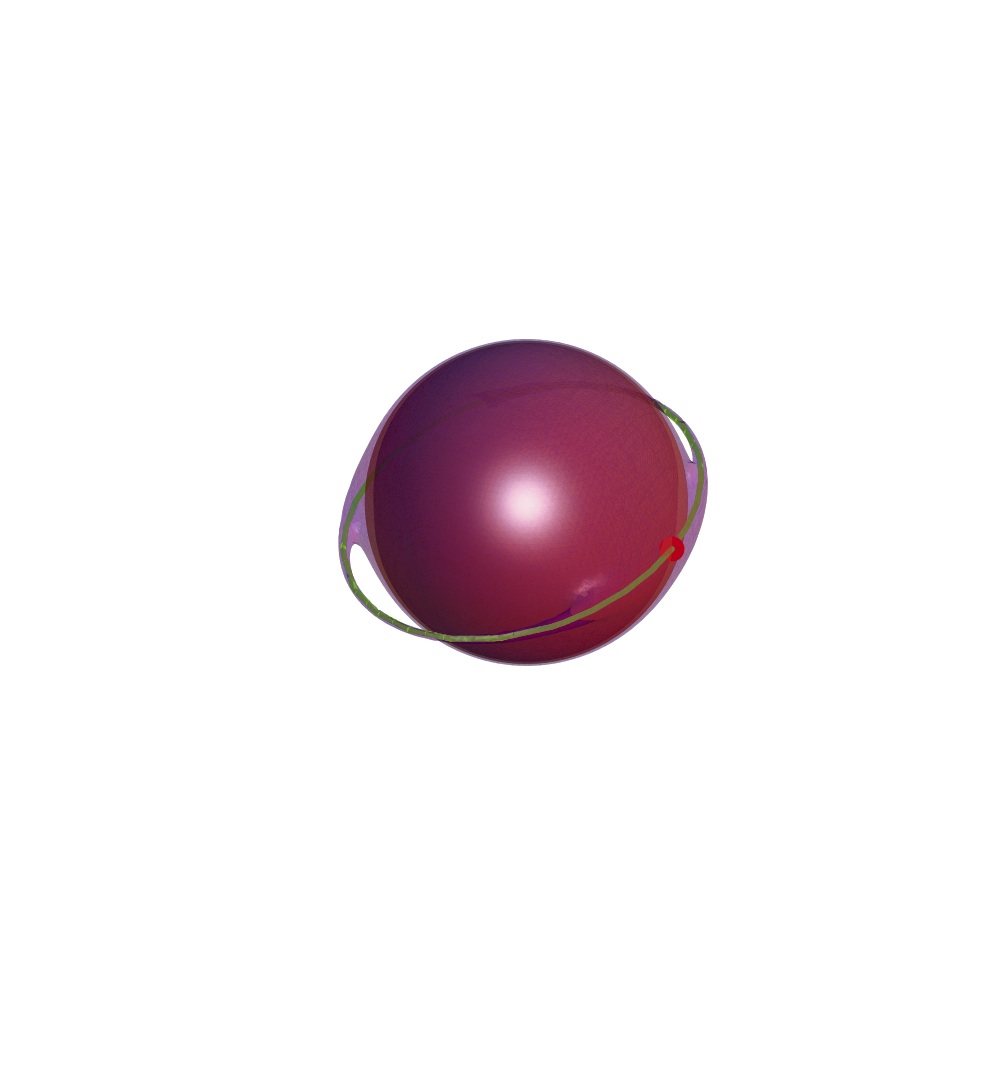}
 \vspace{-5cm}
 \caption{The approximating set $\tilde{V}_{(2)}$.}
  \label{fig:V2}
 \end{figure}
 
 \begin{figure}[h!]
 \includegraphics[scale=1.0]{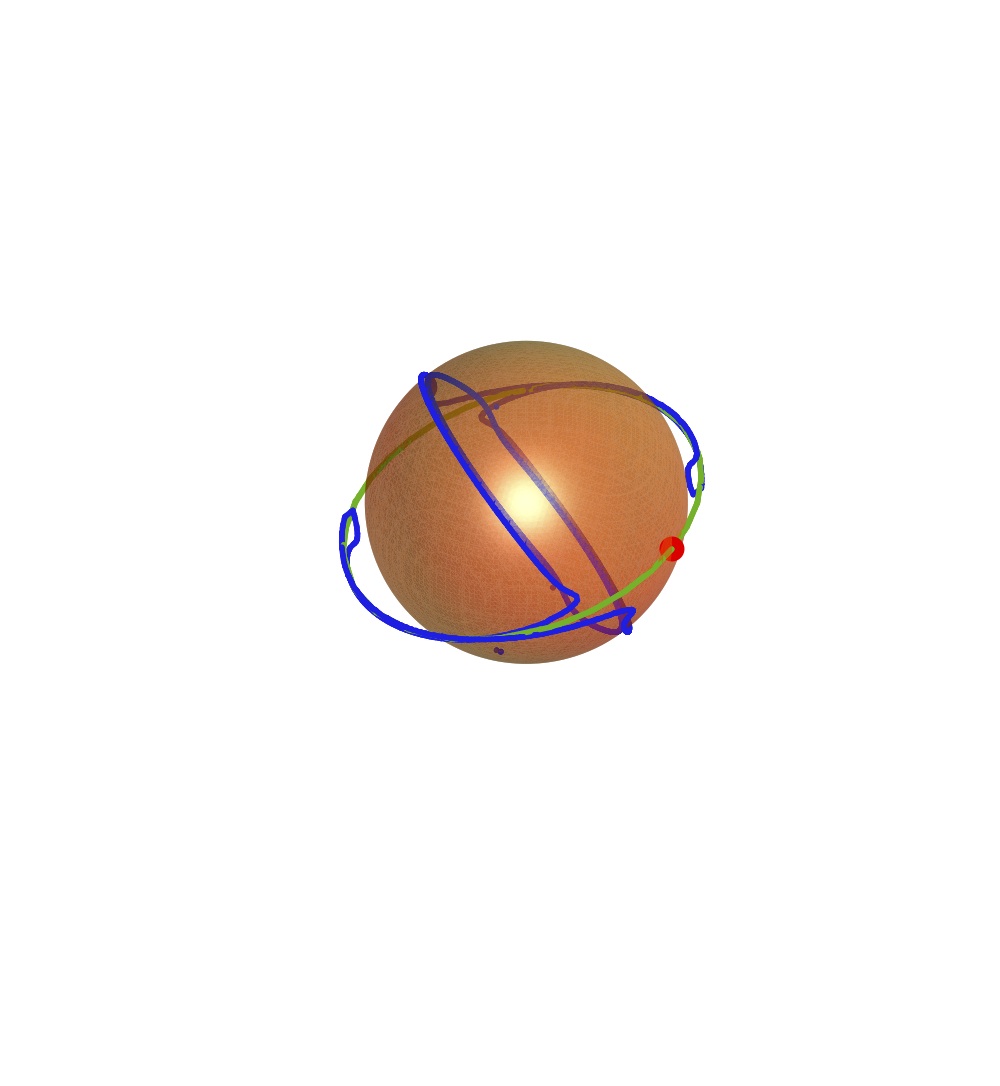}
 \vspace{-5cm}
 \caption{The approximating set $\tilde{V}_{(1)}$ in blue.}
 \label{fig:V1}
 \end{figure}
 
 \begin{figure}[h!]
 \includegraphics[scale=1.0]{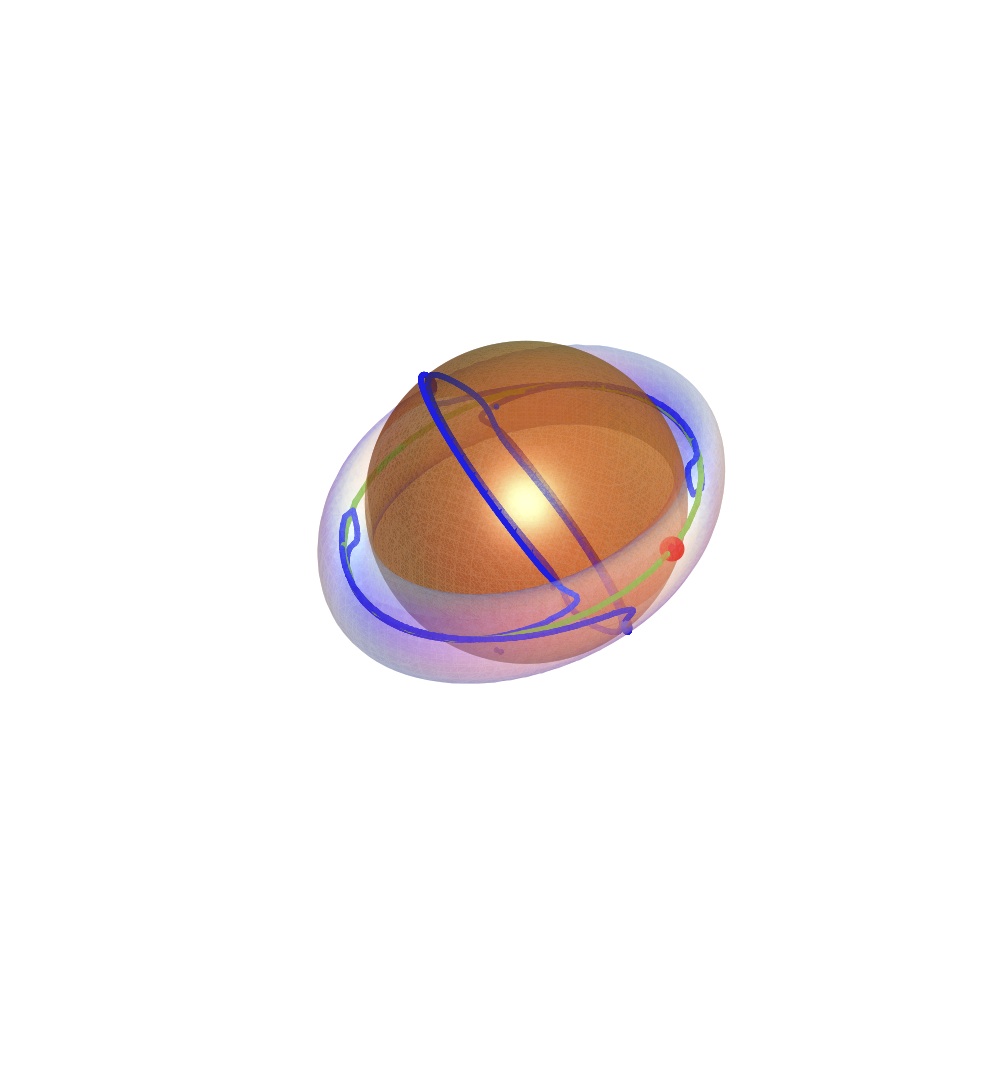}
 \vspace{-5cm}
 \caption{The approximating set $\tilde{V}_{(11)}$.}
  \label{fig:V11}
 \end{figure}
 \end{example}
 
 \begin{figure}[h!]
 \includegraphics[scale=1.0]{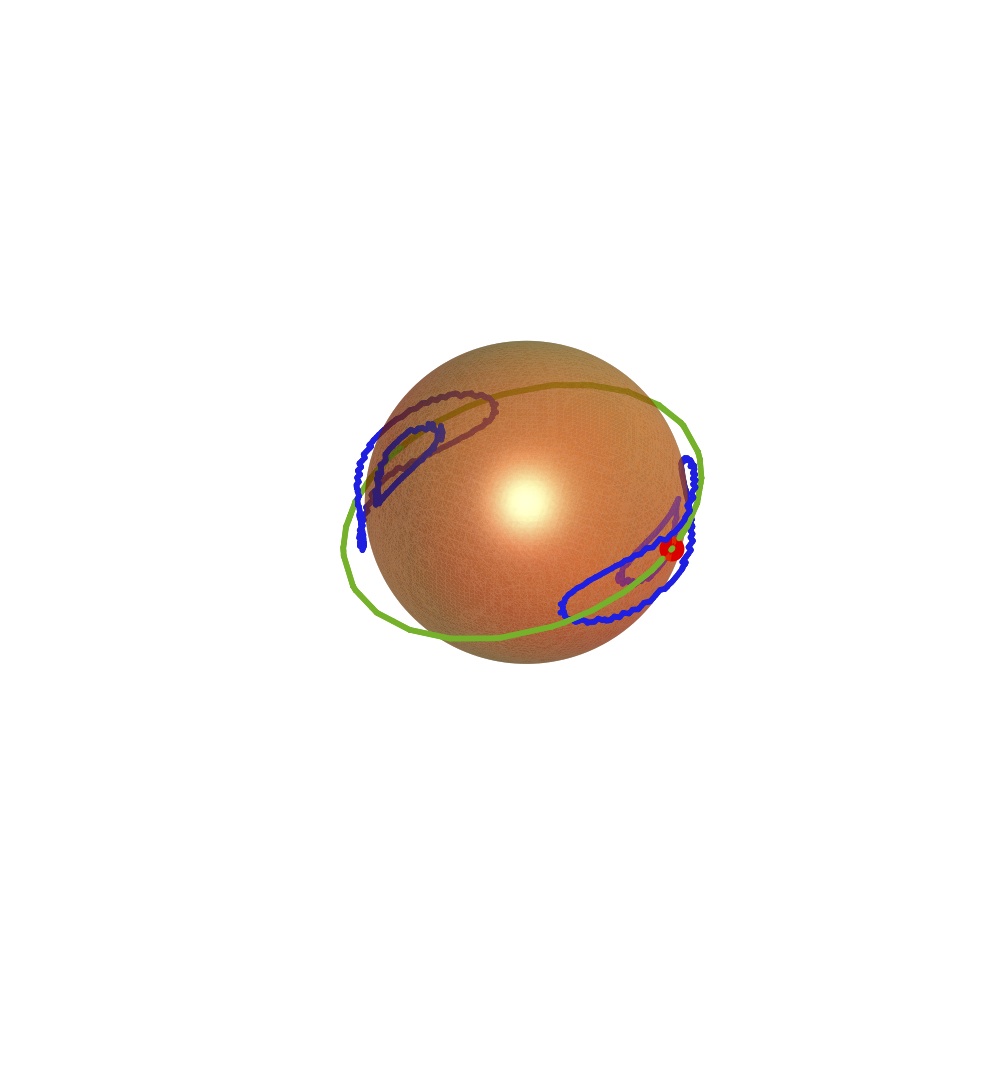}
 \vspace{-5cm}
 \caption{The approximating set $\tilde{V}_{(21)}$ in blue.}
 \label{fig:V21}
 \end{figure}
  
In Figure \ref{fig:V2} we depict the approximating set $\tilde{V}_{(2)}$. Notice that in this example 
 $V_{(2)} \subset \lim_{\delta_1} \tilde{V}_{(2)} \subset V_1$. The first inclusion is
 proper, while the second one is an equality.
 
 In Figure \ref{fig:V1} we depict  the approximating set $\tilde{V}_{(1)}$ (in blue). Note that 
 $V_{(1)} \subset \lim_{\delta_1} \tilde{V}_{(1)} \subset V_1$.  Observe that both inclusions are 
 proper in this case. The image of the curve $\tilde{V}_{(1)}$ (shown in blue in  Figure \ref{fig:V1}) under the $\lim_{\delta_1}$ map
 contains the curve $C$ (shown in green), as well as an additional part contained in the sphere $S$.
 
 The set $\tilde{V}_{(1,1)}$ is the intersection of $\tilde{V}_{(1)}$ with the set defined by the inequality $\tilde{P}_2 \leq 0$, which is a tube containing the set $C$, and
 $V_{(1,1)}=C \subset \lim_{\delta_2} \tilde{V}_{(1,1)} \subset V_2$. Both inclusions are equalities in this case. 
 This is depicted in Figure \ref{fig:V11}.
 
 In Figure \ref{fig:V21} we depict  the approximating set $\tilde{V}_{(2,1)}$ (in blue). Note that the set 
 $V_{(2,1)} = \{(\pm 1,0,0)\}$, 
 and we have the inclusions 
 $V_{(21)} \subset \lim_{\delta_2} \tilde{V}_{(21)} \subset V_2$.  Observe that both inclusions are 
 proper in this case.

We now prove Proposition \ref{prop:main}. While reading the proof particular attention should be paid to the
ordering of the infinitesimals, $\delta_i,\eta_i,\zeta_i, 1\leq i \leq j$ which plays a crucial role.

\begin{proof}[Proof of Proposition \ref{prop:main}]
 We first prove
the inclusion $V_{\tau} \subset W_{\tau}$.

Let $x \in V_{\tau}$ with $\dim^{(j)} (x) = \tau$.  We will prove that $x \in
W_{\tau}$ which suffices to prove the inclusion $V_{\tau} \subset W_{\tau}$,
since $W_{\tau}$ is closed and $V_{\tau}$ is the closure of the set of points
$y$ with $\dim^{(j)} (y) = \tau$. The proof of the claim that $x \in W_{\tau}$
is by induction on $j$. Suppose the claim holds for $j - 1$. There are two
cases to consider.
\begin{enumerate}
  \item\label{item:main1}
   $\tau_j = \tau_{j - 1}$: The induction hypothesis implies that $x
  \in \lim_{\delta_{j - 1}} \tilde{V}_{\sigma'}$, where $\sigma' \in A_{j -
  1}$ with $\sigma'_{j - 1} = \tau_{j - 1} = \tau_j$, and $\sigma'_i \leq
  \tau_i$ for $1 \leq i < j - 1$. Let $\alpha \in I_{j - 1} (\sigma')$ be such
  that $x \in \lim_{\delta_{j - 1}} (V_{\sigma', j - 1}^{\alpha})$. Hence,
  there exists $x' \in V_{\sigma', j - 1}^{\alpha}$ such that $\lim_{\delta_{j
  - 1}} x' = x$. Moreover, since, $\bar{P}_j (x) = 0$, we have that
  $\lim_{\delta_{j - 1}} \bar{P}_j (x') = 0$. From the definition of
  $\tilde{P}_j$ and the fact that $\delta_j \gg \delta_{j - 1} > 0$, we obtain
  that $\tilde{P}_j (x') \leq 0$, and hence $x' \in V_{\sigma, j}^{\beta}$,
  and $x \in \lim_{\delta_j} V_{\sigma, j}^{\beta}$ where $\sigma = (\sigma',
  \tau_j)$, and $\beta = (\alpha, - 1)$.
  
  \item \label{item:main2}
   $q = \tau_j < \tau_{j - 1}$: We prove that every neighborhood, $U$, of
  $x$ in $V_j$ contains a point of $W_{\tau}$. Let $U$ be a small enough
  neighborhood of $x$ in $V_j$. Then there exists a non-empty open subset $U'
  \subset U$ such that each $x' \in U'$ is a regular point of $V_j$ of
  dimension $q$. 
  
  For each $x' \in U'$, 
  shrinking $U'$ further if necessary, we have the inequalities
  $q \leq \dim_{x'} V_{j - 1}
  \leq \dim_x V_{j - 1} = \tau_{j - 1}$, the second inequality coming from
  upper semi-continuity property of the 
  dimension function.
  There are two subcases.
  
  \noindent{{\bf Case (a)}}
  If there
  exists $x' \in U'$, with $\dim_{x'} V_{j - 1} = q = \dim_{x'} V_j$, we
  are reduced to Case (\ref{item:main1}) as follows. Let $\sigma = (\tau_1,\ldots,\tau_{j-2},q,q)$.
  Then, $\sigma \leq \tau$, and using Case (\ref{item:main1}), 
  $x' \in W_\sigma$, and $W_\sigma \subset W_\tau$ (using Lemma \ref{lem:main2}).
  
  \noindent{{\bf Case (b)}}
  We assume that
  $q < \dim_{x'} V_{j - 1} \leq
  \tau_{j - 1}$ for each $x' \in U'$. Using the genericity of the given
  co-ordinates and shrinking $U'$ if necessary by subtracting a Zariski closed
  set of co-dimension at least one we can assume that the tangent space
  $T_{x'} V_j$ is transversal to $\pi^{- 1}_{[1, q]} (z')$ (recall Notation
  \ref{not:projection}), where $z' = \pi_{[1. q]} (x')$, and hence in
  particular that $x'$ is an isolated point of $(V_j)_{z'}$ for all $x' \in U'$.
  
  Shrinking $U'$ further if necessary we can also assume that $x'$ is not an isolated point
  of $(V_{j - 1})_{z'}$ where $z' = \pi_{[1. q]} (x')$ for all $x' \in U'$. To
  see this suppose that there exists a non-empty open subset $U''$ of $U'$
  such that for all $x'' \in U''$, $x''$ is an isolated point of  $(V_{j -
  1})_{z''}$ where $z'' = \pi_{[1. q]} (x'')$. Then, there exists for any $x''
  \in U''$ an open neighborhood $W$ of $x''$ in $V_{j - 1}$ contained in
  $(V_{j - 1})_{\pi_{[1, q] (U'')}}$ such that the dimension of $W$ is $\leq
  q$, which is contrary to our assumption. 
  
  Now for each $x' \in U'$, since $x'$ is an isolated point of $(V_j)_{z'}$,
  and $\lim_{\delta_j}  ((\Bas(\emptyset, \{\tilde{P}_j\}))_{z'})_b = (V_j)_{z'}$,
  there exists a 
  unique
  semi-algebraically connected component of $((\Bas(\emptyset,
  \{\tilde{P}_j\}))_{z'})_b$, and hence  of $\tmop{Zer} (\tilde{P}_j,
  \R \langle \delta_j \rangle^k)_{z'}$ (which we will denote by  $Z(x')$)
  such that $\lim_{\delta_j} Z(x') = x'$. Since $x'$ is not an isolated point
  of $(V_{j - 1})_{z'}$, 
  $\tmop{Ext} ((V_{j - 1})_{z'},
  \R \langle \delta_j \rangle) \cap Z(x') \neq
  \emptyset$.

  We claim that there exist, 
  $\sigma'  = (\sigma'_1,\ldots,\sigma'_{j-1}) \in A_{j - 1}$,  
  $\sigma' \leq \tau' :=(\tau_1, \ldots, \tau_{j -1})$,
 $\alpha \in I_{j - 1} (\sigma')$, 
 $x' \in U'$, $z' = \pi_{[1. q]} (x')$, 
 such that 
 \begin{equation*}
 \lim_{\delta_{j-1}} ((V_{\sigma', j - 1}^{\alpha})_{z'})_b
 \cap  
 \tmop{Ext}((V_{\tau'})_{z'}, \R \langle \delta_j \rangle)\cap Z(x')
 \neq \emptyset,
 \end{equation*}
 where $Z(x')$ is the unique semi-algebraically connected component of 
 \[
 \tmop{Zer} (\tilde{P}_j,
  \R \langle \delta_j \rangle^k)_{z'}
 \]
  such that $\lim_{\delta_j} Z(x')= x'$ (see previous paragraph).
  
  To see this let $\sigma'$ be a minimal element in $A_{j-1}$ such that 
  $U'' := V_{j-1}^{\leq \sigma'} \cap U' \neq \emptyset$. By Proposition \ref{prop:uppersc},  
  (see also Example \ref{eg:main} and  Remark \ref{rem:ex-main} following it)
  $U''$ is a non-empty open subset of $U'$. 
  Now suppose that  
  \[
  \lim_{\delta_{j -1}} ((V_{\sigma', j - 1}^{\alpha})_{z'})_b \cap \tmop{Ext}
  ((V_{\sigma'})_{z'}, \R \langle \delta_j \rangle)
  \cap Z(x') = \emptyset,
  \] 
  for every $x' \in U''$ and $z' = \pi_{[1,q]}(x')$. 
  Now,
  $\pi_{[1,q]}(U')$, and hence  $\pi_{[1,q]}(U'')$, is a non-empty open subset of $\R^{[1,q]}$,
  since the map $\pi_{[1,q]}|{U'}$ is a semi-algebraic diffeomorphism by the 
  semi-algebraic implicit function theorem, and
  the fact that $T_{x'} V_j$ is transversal to $\pi^{- 1}_{[1, q]} (\pi_{[1,q]}(x'))$ for every 
  $x' \in U'$ (see above).
  This contradicts the inductive
  hypothesis,  which implies that $V_{\sigma'} \subset W_{\sigma'}$.
  
  Now fix $x', z', \sigma', \alpha$ as above. Notice that since $\lim_{\delta_j} Z(x') = x'$, there exists $r > 0$, such that
  \[
  Z(x') = \tmop{Zer} (\tilde{P}_j, \R \langle
  \delta_j \rangle^k)_{z'} \cap \overline{B_k (x', r)}_{z'}.
  \]
  It now follows from Lemma \ref{lem:technical} 
  (applied after taking $\delta = \delta_j$ and $\bar\eps=(\delta_{j-1},\ldots,\delta_1,\eta_1,\zeta_1,\ldots,\eta_{j-1},\zeta_{j-1})$)
  that 
  \[
  \lim_{\delta_{j - 1}}
  ((V_{\sigma', j - 1}^{\alpha})_{z'})_b \cap Z(x') \neq \emptyset
  \]
  implies that
  \begin{equation}
  \label{eqn:application-technical}
  \tmop{Ext} ((V_{\sigma', j - 1}^{\alpha})_{z'},
  \R') \cap \tmop{Ext}(Z(x'),\R')
  \neq\emptyset,
  \end{equation}
  where
  $\R' = \R\la\delta_j,\ldots,\delta_1,\eta_1,\zeta_1,\ldots,\eta_{j-1},\zeta_{j-1}\ra$.
  Moreover, it is clear that \eqref{eqn:application-technical} implies that
  \[
  \tmop{Ext} ((V_{\sigma', j - 1}^{\alpha})_{z'},
  \R_j) \cap \tmop{Ext}(Z(x'),\R_j)
  \neq\emptyset.
  \]
  Note that the order $\delta_j \gg \delta_{j - 1}$ is important here (cf. Remark
  \ref{rem:ordering-of-infinitesimals}).
  
  It follows that there exists a semi-algebraically connected component
  $C$ of 
  $
  \tmop{Zer} ((\mathcal{P}^{\alpha}_{\sigma', j - 1}, \tilde{P}_j),
  \R_j^k)_{z'},
  $
  such that $x' \in \lim_{\delta_j} C$, and $C \subset (V_{\sigma', j - 1}^{\alpha})_{z'}$. 
  Moreover, using the fact that $z' \in \R^{[1, q]}$, and applying
  Proposition \ref{prop:transversal} 
  with $\eta=\eta_j,\delta=\delta_j$, and $\bar\eps=(\delta_{j-1},\ldots,\delta_1,\eta_1,\zeta_1,\ldots,\eta_{j-1},\zeta_{j-1})$,
 we deduce that the polynomials in
  $(\mathcal{P}^{\alpha}_{\sigma', j - 1} (z', \cdot), \tilde{P}_j (z',
  \cdot))$ define a non-singular complete intersection of dimension $p - q -
  1$ in $\R_j^{[q + 1, k]}$, where $p = \sigma'_{j
  - 1}$. Let $\tilde{\mathcal{P}} =\mathcal{P}^{\alpha}_{\sigma', j - 1}$, and
  let $\mathcal{F}$ be the tuple of polynomials defined by
  \eqref{eqn:FF}. Then, there exists a semi-algebraically connected component
  $\tilde{C}$ of $\tmop{Zer} (\mathcal{F},
  \R_j^k)_{z'}$ such that $x' \in \lim_{\delta_j}
  \tilde{C}$. There are a finite number of $X_k$-critical points (all of which
  are simple) on $\tilde{C}_z$ by Remark \ref{rem:generic} and Proposition
  \ref{prop:transversal}. If $(z', w') $, $w' \in
  \R_j^{[q + 1 +, k]}$, is one such critical
  point, then $(z', w')$ is contained in the finite constructible
  set $C_J (\mathcal{F})$ (cf. \eqref{eqn:CJ} and part \ref{item:prop:transversal3}. of Proposition \ref{prop:transversal}) for some $J \in \binom{[q + 1, k]}{k - p + 1}$, and
  such that $w'$ is a simple zero of the system $ \mathcal{\mathcal{F}}_J
  (z', \cdot)$. Hence,  applying Proposition \ref{prop:simple2} 
  (with the field of coefficients $\R$ in the Proposition \ref{prop:simple2}  taken to be the real closed field
  $\R\la\delta_j,\ldots,\delta_1,\eta_1,\zeta_1,\eta_{j-1},\zeta_{j-1},\eta_j\ra$ and $\zeta=\zeta_j$) we have that
  there exists a simple
  zero, $w''$, of the system $\mathcal{P}^{\beta}_{\sigma, j} (z, \cdot)$ (cf.
  (\ref{eqn:deformedF})) where $\sigma = (\sigma', \tau_j)$ and $\beta =
  (\alpha, J)$, such that $\lim_{\zeta_j} w'' = w'$. Clearly, then $x'' =
  (z', w'') \in V_{\sigma, j}^{\beta}$, , and $x' = \lim_{\delta_j} x''$ and
  thus $x' \in \lim_{\delta_j} V_{\sigma, j}^{\beta}$. Notice that $\sigma_j =
  \tau_j$ and $\sigma \leq \tau$.
\end{enumerate}

The inclusion $\lim_{\delta_j} \tilde{V_{\tau}} \subset V_j$, from which
the second inclusion $W_{\tau} \subset V_j$ follows immediately, is due to the
fact that for each $\beta \in I_j (\tau)$, $V_{\tau, j}^{\beta}$ is either
contained in the part of the semi-algebraic set defined by $\tilde{P}_j
\leq 0$ which is bounded over $\R$, or in the
algebraic variety $\tmop{Zer} (\tilde{P}_j,
\R_j^k)_b$ depending on whether $\tau_{j - 1} =
\tau_j$ or $\tau_{j - 1} > \tau_j$ respectively. It is clear from definition
of $\tilde{P}_j$, that the images under $\lim_{\delta_j}$  of the last two sets
are contained in $V_j$.
\end{proof}

The following slight refinement of Proposition \ref{prop:main} is required to
ensure that the degree of the last polynomial does not enter the bound with a
factor of $(k - \tau_{i - 1} - 1)$ as is the case of the other degrees $d_{i
}$, with $i < \ell$, but rather just as $d_{\ell}$. This slight improvement
is possible since we do not need to ensure that the dimension of the
approximating varieties drops appropriately (to $k_{\ell}$) when we
approximate the last variety $V_{\ell}$. If we were not interested in
obtaining the tightest possible dependence on $k$ in the multiplicative factor
in the bound (the factor that is independent of the degrees), then this
refinement would not have been necessary. However, in order to ensure that the
results in the current paper properly generalize the results in
{\cite{Barone-Basu11a}} we need to take this extra care.

\begin{notation}
  For all $\sigma = (\sigma_1, \ldots, \sigma_j) \in A_j, 2\leq j\leq \ell$, denote by 
  \[
  \hat{\sigma}
  = (\sigma_1, \ldots, \sigma_{j - 1}, \sigma_{j - 1}).
  \]
\end{notation}

\begin{corollary}
  \label{cor:main} For all $\tau \in A_j$,
  \[ V_{\tau} \subset W_{\tau}' \subset V_j , \]
  where
  \[ 
  W_{\tau}' = \bigcup_{\substack {\sigma=(\sigma_1,\ldots,\sigma_j) \in A_j\\ \sigma_i \leq \tau_i,1 \leq i < j, \sigma_j = \tau_j}} 
  \lim_{\delta_j} \tilde{V}_{\hat{\sigma}}.
  \]
\end{corollary}

\begin{proof} 
It is clear from the definition that for all
$\sigma \in A_j$
\[ \tilde{V}_{\sigma} \subset \tilde{V}_{\hat{\sigma}} , \]
and that $\lim_{\delta_j} \tilde{V}_{\hat{\sigma}} \subset V_j$. The corollary now
follows from Proposition \ref{prop:main}.
\end{proof}

\begin{corollary}
  \label{cor:main2}
 \[ b_0 (V_{\ell}) \leq \sum_{\tau \in A_{\ell}} \sum_{\beta \in I_{\ell}
     (\hat{\tau})} b_0 (V_{\hat{\tau}, \ell}^{\beta}) . \]
\end{corollary}

\begin{proof}
Follows immediately from Corollary \ref{cor:main}
after noting that (using Proposition \ref{prop:cover})
\[ V_{\ell} = \bigcup_{\tau \in A} V_{\tau} . \]
\end{proof}

Following notation introduced above we have the following proposition.

\begin{proposition}
  \label{prop:quantitative} Let $\tau \in A_j $, $\tau_{j - 1} = p$, and
  $\alpha \in I_j (\hat{\tau})$. 
  \begin{enumerate}[1.]
    \item\label{item:prop:quantitative1} Then $\tmop{card} \mathcal{P}^{\alpha}_{\hat{\tau}, j} = k - p$.
    \item \label{item:prop:quantitative2} Suppose that $\mathcal{P}^{\alpha}_{\hat{\tau}, j} = (P_1, \ldots, P_{k -
    p})$. Let for $1 \leq i \leq j - 1$, $\ell_i = \tau_{i - 1} - \tau_i$,
    with the convention that $\tau_0$=k, and $L_i = \sum_{h = 1}^i \ell_h$.
    Then for each $i , 1 \leq i < j$, the degrees of the polynomials
    $P_{L_{i - 1} + 1}, \ldots, P_{L_i}$ are bounded by $(k - \tau_{i - 1} +
    1) d_i \leq (k + 1) d_i$.
    
    \item\label{item:prop:quantitative3} $\deg({Q}^{\alpha}_{\hat{\tau},
    j})\leq d_{\ell}$.
    
    \item \label{item:prop:quantitative4} $b_0 (V^{\alpha}_{\hat{\tau}, j}) \leq 
    b_0 (\tmop{Zer}(\mathcal{P}^{\alpha}_{\hat{\tau}, j}\cup\{ {Q}^{\alpha}_{\tau, j}\},\R_j^k)_b)+
    b_0 (\tmop{Zer}(\mathcal{P}^{\alpha}_{\hat{\tau}, j},\R_j^k)_b)$.
  \end{enumerate}
\end{proposition}

\begin{proof}
Follows from the definitions of the tuples
$\mathcal{P}^{\alpha}_{\hat{\tau}, j}$, and the polynomials 
${Q}^{\alpha}_{\hat{\tau}, j}$ (see Definition \ref{def:tuple}), as well as 
Proposition \ref{prop:algebraic}.
\end{proof}

\subsection{Bounds on the 
$0$-th Betti number
of  non-singular complete intersections}
\label{subsec:bounds-on-Betti}
The following proposition appears in \cite{Barone-Basu11a}, and is a
consequence of the classical formula for the Euler-Poincar\'e characteristic of
non-singular complex projective intersections and the Smith inequality.

\begin{proposition}
  \label{prop:hirzebruch} Let $\mathcal{F} = \{ F_1, \ldots, F_m \} \subset
  \R [X_1, \ldots, X_k]$ with $\deg (F_i) =
  d_i$, $d_1 \leq d_2 \leq \cdots \leq d_m$. Moreover, assume that
  $\mathcal{F}^h = \{ F_{1}^h, \ldots, F_{m}^h \}$ defines a
  non-singular complete intersection in
  $\mathbb{P}^k_{\C}$. Then,
  \[ b_0 (\ensuremath{\operatorname{Zer}} (\mathcal{F},
     \R^k)_b) \leq \binom{k + 1}{m + 1} d_1
     \cdots d_{m - 1} d_m^{k - m + 1} + 2 (k - m + 1) . \]
\end{proposition}

\begin{remark}
  \label{rem:Katz} We note that in Proposition \ref{prop:hirzebruch} if the
  polynomials in $\mathcal{F}$ do not define a non-singular complete
  intersection, it is still possible to bound the sum of the Betti numbers of
  the corresponding complex variety by $O (1)^m O (m d_m)^k$ using a result of
  Katz {\cite{Katz2001}}, which in turn uses previous results of Bombieri
  {\cite{Bombieri78}}, and Adolphson and Sperber {\cite{Adolphson-Sperber}}.
  These results use the theory of exponential sums over finite fields, and are
  of a much deeper nature than the classical formula giving the Betti numbers
  in terms of the degree sequence in the non-singular complete intersection
  case which is used to prove Proposition \ref{prop:hirzebruch}. However, the
  results of Katz {\cite{Katz2001}} which do not assume non-singularity and
  are very general, do not have the finer dependence on the degree sequence
  (see the bound given above), and this finer dependence on the degree
  sequence is the key point in Proposition \ref{prop:hirzebruch} above.
\end{remark}

\begin{corollary}
  \label{cor:hirzebruch} For each $\tau = (\tau_1, \ldots, \tau_{\ell}) \in
  A_{\ell}$ and $\alpha \in I_{\ell} (\hat{\tau})$ and $\mathcal{Q} \subset
  \{{Q}^{\alpha}_{\hat{\tau}, \ell}\}$,
  \[ b_0 (\ensuremath{\operatorname{Zer}} ((\mathcal{P}^{\alpha}_{\hat{\tau},
     \ell}\cup \mathcal{Q})_b, \R_{\ell}^k)) \leq O
     (1)^k d_{\ell}^{\tau_{\ell - 1}} \prod_{1 \leq i < \ell} ((k - \tau_{i -
     1} + 1) d_i)^{\tau_{i - 1} - \tau_i} . \]
\end{corollary}

\begin{proof}
Follows from parts \ref{item:prop:quantitative1}.,  \ref{item:prop:quantitative2}. and \ref{item:prop:quantitative3}. of Proposition
\ref{prop:quantitative}, and Proposition \ref{prop:hirzebruch}.
\end{proof}

It now follows from Corollary \ref{cor:hirzebruch} and part \ref{item:prop:quantitative4}.  of Proposition
\ref{prop:quantitative} that

\begin{corollary}
  \label{cor:main3} For each $\tau = (\tau_1, \ldots, \tau_{\ell}) \in
  A_{\ell}$ and $\alpha \in I_{\ell} (\hat{\tau})$
  \begin{eqnarray*}
    b_0 (V^{\alpha}_{\hat{\tau}, \ell}) & \leq & O (1)^k d_{\ell}^{\tau_{\ell - 1}}
    \prod_{1 \leq i < \ell} ((k - \tau_{i - 1} + 1) d_i)^{\tau_{i - 1} -
    \tau_i} .
  \end{eqnarray*}
\end{corollary}

Let $\tau \in A_{\ell}  $ and $d_1, \ldots, d_{\ell}$ satisfy the hypothesis
of Theorem \ref{thm:refined-algebraic}.

\begin{lemma}
  \label{lem:maximized} Then,
  \begin{eqnarray*}
    \frac{d_{\ell}^{\tau_{\ell - 1}} \prod_{1 \leq i < \ell} ((k - \tau_{i -
    1} + 1) d_i)^{\tau_{i - 1} - \tau_i}}{d_{\ell}^{k_{\ell - 1}} \prod_{1
    \leq i < \ell} ((k - k_{i - 1} + 1) d_i)^{k_{i - 1} - k_i}} & \leq & O
    (k)^k .
  \end{eqnarray*}
\end{lemma}

\begin{proof}
 Using the inequality that for $2 \leq i \leq
\ell,$
\begin{eqnarray*}
  \frac{d_{i - 1}}{d_i} & \leq \frac{1}{k + 1} \leq \frac{1}{k - k_{i - 2} +
  1} & 
\end{eqnarray*}
we get that the expression on the left hand side of the proposition is bounded
by
\[ \frac{\prod_{1 \leq i < \ell} (k - \tau_{i - 1} + 1)^{\tau_{i - 1} -
   \tau_i}}{\prod_{1 \leq i < \ell} (k - k_{i - 1} + 1)^{k_{i - 1} - \tau_i}}
   . \]
The sum of the various exponents of the numerator is
\begin{eqnarray*}
  \sum_{i = 1}^{\ell - 1} (\tau_{i - 1} - \tau_i) & = & \tau_0 - \tau_{\ell -
  1} \leq k,
\end{eqnarray*}
and for each $i , 1 \leq i < \ell,$ $(k - \tau_{i - 1} + 1) \leq (k
+ 1)$. The denominator is a non-zero integer.
\end{proof}

We next bound the cardinality of the index set $A_{\ell}$.

\begin{lemma}
  \label{lem:cardinalityofA} The cardinality of $A_{\ell}$ is bounded by
  \[ O (1)^{k + \ell} . \]
\end{lemma}

\begin{proof}
The number of tuples $\tau = (\tau_1, \ldots,
\tau_{\ell})$ in which $k \geq \tau_1 > \tau_2 > \cdots > \tau_{\ell} \geq 0$
is bounded by the volume of the corresponding $\ell$-dimensional simplex in
$\mathbb{R}^{\ell}$ which is equal to $\frac{(k + 1)^{\ell}}{\ell !}$.
Allowing some of the $\tau_i$'s to be equal, the number of tuples is bounded
by
\[ \sum_{0 \leq i \leq \ell} \binom{\ell}{i} \frac{(k + 1)^{\ell - i}}{(\ell
   - i) !} \leq 2^{\ell} \sum_{0 \leq i \leq \ell} \frac{(k + 1)^{\ell -
   i}}{(\ell - i) !} = O (1)^{k + \ell} . \]
\end{proof}

\begin{lemma}
  \label{lem:cardinalityofI} For each $\tau = (\tau_1, \ldots, \tau_{\ell})$
  the cardinality of the index set $I_{\ell} (\tau)$ is bounded by
  \[  \]
  \begin{eqnarray*}
    (k - \tau_{\ell} + 1) \binom{k - \tau_{\ell}}{\tau_1 - \tau_2, \ldots,
    \tau_{\ell - 1} - \tau_{\ell}} . &  & 
  \end{eqnarray*}
  
\end{lemma}

\begin{proof}
 It is clear from the definition that the
cardinality of the index set $I_{\ell} (\tau)$ is bounded by

\begin{eqnarray*}
  \prod_{1 \leq j \leq \ell} \binom{k - \tau_j + 1}{k - \tau_{j - 1} + 1} & =
  & \frac{(k - \tau_{\ell} + 1) !}{(\tau_0 - \tau_1) ! (\tau_1 - \tau_2) !
  \cdots (\tau_{\ell - 1 -} \tau_{\ell}) !}\\
  & = & (k - \tau_{\ell} + 1) \binom{k - \tau_{\ell}}{\tau_0 - \tau_1, \tau_1
  - \tau_2, \ldots, \tau_{\ell - 1} - \tau_{\ell}} .
\end{eqnarray*}
\end{proof}

\subsection{Proof of Theorem \ref{thm:refined-algebraic}}
\label{subsec:proof-of-theorem-algebraic}
We now prove Theorem \ref{thm:refined-algebraic}.
\label{subsec:proof-of-theorem-algebraic}
\begin{proof}[Proof of Theorem \ref{thm:refined-algebraic}]
We first
prove the theorem in case $V_0$ is bounded. It follows from Corollary
\ref{cor:main2} and Corollary \ref{cor:main3} that
\begin{eqnarray*}
  b_0 (V_{\ell}) & \leq & \sum_{\tau \in A_{\ell}}  \sum_{\alpha \in I_{\ell}
  (\hat{\tau})} \left( O (1)^k d_{\ell}^{\tau_{\ell - 1}} \prod_{1 \leq i < \ell}
  ((k - \tau_{i - 1} + 1) d_i)^{\tau_{i - 1} - \tau_i} \right) .
\end{eqnarray*}
Using Lemma \ref{lem:cardinalityofI} to bound the cardinality of the index
set $I_{\ell} (\hat{\tau})$, we get that the right hand side of the above
inequality is bounded by
\[ O (1)^k \sum_{\tau \in A_{\ell}} F (k, \tau) \left( d_{\ell}^{\tau_{\ell -
   1}} \prod_{1 \leq i < \ell} ((k - \tau_{i - 1} + 1) d_i)^{\tau_{i - 1} -
   \tau_i} \right), \]
where
\[ F (k, \tau) = (k - \tau_{\ell - 1} + 1) \binom{k - \tau_{\ell - 1}}{\tau_0
   - \tau_1, \tau_1 - \tau_2, \ldots, \tau_{\ell - 2} - \tau_{\ell - 1}} . \]
The theorem in the bounded case now follows from Lemma \ref{lem:maximized} and
Lemma \ref{lem:cardinalityofA}.

In the general case, we first replace the given sequence of polynomials $Q_1,
\ldots, Q_{\ell}$, by a new sequence, $Q_0, Q_1, \ldots, Q_{\ell}$, where
\[ Q_0 = \sum_{i = 1}^{k + 1} X_i^2 - \Omega, \]
where $\Omega$ is infinitely large and positive over
$\R$. For each $i, 0 \leq i \leq \ell$, defining
$\hat{\mathcal{Q}}_i = \{ Q_0, \ldots, Q_i \}$, and $\hat{V}_i = \tmop{Zer}
(\hat{\mathcal{Q}}_i, \R \langle 1 / \Omega
\rangle^{k + 1})$ we have that each $\widehat{V_i}$ is bounded over
$\R \langle 1 / \Omega \rangle$, and also that \
$b_0 (V_{\ell}) \leq b_0 (\hat{V}_{\ell})$. Applying the same arguments as
in the bounded case we obtain that
\[ b_0 (\hat{V}_{\ell}) \leq O (1)^k \sum_{\tau = (\tau_{- 1}, \tau_0,
   \ldots, \tau_{\ell - 1})} F (k, \tau) \left( d_{\ell}^{\tau_{\ell - 1}}
   \prod_{1 \leq i < \ell} ((k - \tau_{i - 1} + 1) d_i)^{\tau_{i - 1} -
   \tau_i} \right), \]
where the sum is taken over all $\tau \in \mathbb{N}^{\ell}$, with $k + 1 =
\tau_{- 1} > k = \tau_0 \geq \tau_1 \cdots \geq \tau_{_{\ell - 1}} \geq 0$,
and $\tau_i \leq k_i$, for each $i, 1 \leq i < \ell$, and
\[ F (k, \tau) = (k - \tau_{\ell - 1} + 1) \binom{k - \tau_{\ell - 1}}{\tau_0
   - \tau_1, \tau_1 - \tau_2, \ldots, \tau_{\ell - 2} - \tau_{\ell - 1}} . \]
Notice that since the local dimension of the variety $\hat{V}_0$ is
constant, it suffices to fix $\tau_0 = k$ in the sum above, and the
contribution of the degree of the polynomial $Q_0$ (note that $\deg(Q_0) =2$)  gets absorbed into the $O
(1)^k$ term.

\end{proof}

\subsection{Proof of Theorem \ref{thm:refined-semi-algebraic}}
\label{subsec:proof-of-theorem-semi-algebraic}
We now prove Theorem \ref{thm:refined-semi-algebraic}.

We introduce a new family of polynomials defined as follows:
\[ \tilde{\mathcal{P}} = \bigcup_{1 \leq i \leq s} \{ P_i \pm \varepsilon
   \gamma_i, P_i \pm \delta \gamma_i \}, \]
where $\varepsilon, \delta, \gamma_1, \ldots, \gamma_s$ new variables. 

For any
subset $I = \{ (\epsilon_1, \sigma_1, i_1), \ldots, (\epsilon_{m,} \sigma_m,
i_m) \} \subset \{ + 1, - 1 \} \times \{ \varepsilon, \delta \} \times \{ 1,
\ldots s \}$,
we denote by $\tilde{\mathcal{P}}_I$ the subset of
$\tilde{\mathcal{P}}$ defined by
\begin{eqnarray*}
  \tilde{\mathcal{P}}_I & = & \bigcup_{1 \leq j \leq m} \{ P_{i_j} +
  \epsilon_j \sigma_j \gamma_{i_j} \} .
\end{eqnarray*}
Let $\R'$ denote the real closed field
$\R \langle \varepsilon, \delta, \gamma_1,
\ldots \gamma_s \rangle$.

\begin{proposition}
  \label{prop:genpos} For each $I \subset \{ + 1, - 1 \} \times \{
  \varepsilon, \delta \} \times \{ 1, \ldots, s \}$, the dimension of the
  variety $\ensuremath{\operatorname{Zer}} (\tilde{\mathcal{P}},
  \R'^k) \cap \ensuremath{\operatorname{Ext}}
  (V_{\ell}, \R')$ is at most $k_{\ell}
  -\ensuremath{\operatorname{card}}I$. In particular,
  $\ensuremath{\operatorname{Zer}} (\tilde{\mathcal{P}},
  \R'^k) \cap \ensuremath{\operatorname{Ext}}
  (V_{\ell}, \R')$ is empty if
  $\ensuremath{\operatorname{card}}I > k_{\ell}$.
\end{proposition}

\begin{proof} It follows immediately from the fact that the
various $\gamma_i$'s are algebraically independent over
$\R$.
\end{proof}

\begin{notation}
\label{not:weak}
For any finite family $\mathcal{F} \subset \R
[X_1, \ldots, X_k]$ we call a formula 
\[
\bigwedge_{F \in \mathcal{F}} F
\sigma_F 0,
\] 
where each $\sigma_F \in \{ \geq, \leq \}$,  
a \emph{weak sign condition} on $\mathcal{F}$. 
\end{notation}

\begin{proposition}
  \label{prop:weak} Let $V_1$ be bounded, and let $\sigma \in \{ - 1, 0, 1
  \}^{\mathcal{P}}$ and $C$ a semi-algebraically connected component of
  $\ensuremath{\operatorname{Reali}} (\sigma, V_{\ell}) \subset
  \R^k$. Then, there exists a weak sign
  condition $\tilde{\sigma}$ on $\tilde{\mathcal{P}}$,  and a semi-algebraically
  connected component $\tilde{C}$ of 
  \[
  \ensuremath{\operatorname{Reali}}
  (\tilde{\sigma}, \ensuremath{\operatorname{Ext}} (V_{\ell},
  \R'))
  \] 
  such that 
  \[
  \lim_{\delta} \tilde{C}
  \subset \ensuremath{\operatorname{Ext}} (C,
  \R').
  \]
  \end{proposition}

\begin{proof}
 The proof is similar to the proof of Proposition 4 in
{\cite{BPR95a}} and omitted.
\end{proof}

The following proposition occurs in {\cite{BPRbook2}} (Proposition 13.1).

\begin{proposition}
  \label{prop:semi-algebraic} Let $V_1$ be bounded and let $\mathcal{F}
  \subset \R [X_1, \ldots, X_k]$ be a finite set
  of polynomials and $\tilde{\sigma}$ a weak sign condition on $\mathcal{F}$.
  Let $C$ be a semi-algebraically connected component of
  $\ensuremath{\operatorname{Reali}} (\tilde{\sigma},
  \ensuremath{\operatorname{Ext}} (V_{\ell}, \R'
  ))$. Then there exists a subset $\mathcal{F}' \subset \mathcal{F}$, and a
  semi-algebraically connected component $D$ of
  $\ensuremath{\operatorname{Zer}} (\mathcal{F}',
  \ensuremath{\operatorname{Ext}} (V_{\ell}, R'))$, such that $D \subset C.$
\end{proposition}

\begin{proof}[Proof of Theorem \ref{thm:refined-semi-algebraic}]
In
the case $V_1$ is bounded, using successively Propositions \ref{prop:weak} and
\ref{prop:semi-algebraic} it suffices to bound the total number of
semi-algebraically connected components of the real algebraic sets
\[ \tmop{Zer} (\mathcal{Q}_{\ell} \cup \tilde{\mathcal{P}}_I,
   \R'^k) \]
for subsets $I \subset \{ + 1, - 1 \} \times \{ \varepsilon, \delta \} \times
\{ 1, \ldots s \}$. Moreover, using Proposition \ref{prop:genpos},  the set of
different subsets $I$ that we need to consider is bounded by
\[ \sum_{j = 0}^{k_{\ell}} 4^j \binom{s}{j} = (O (s))^{k_{\ell}} . \]

Notice that each $\tmop{Zer} (\mathcal{Q}_{\ell} \cup \tilde{\mathcal{P}}_I,
\R'^k) = \tmop{Zer} ((Q_1, \ldots, Q_{\ell}, P_I),
\R'^k)$, where $P_I = \sum_{P \in
\tilde{\mathcal{P}}_I} P^2$. 
Also, notice  that
\[
(\deg(Q_1),\ldots, \deg(Q_{\ell}), \deg(P_I)) = (d_1, \ldots, d_{\ell}, 2 d).
\] Now apply Theorem
\ref{thm:refined-algebraic} to finish the proof. In the general case, use the
same technique as in the proof of Theorem \ref{thm:refined-algebraic} to
reduce to the bounded case.
\end{proof}

\begin{proof}[Proof of Theorem \ref{thm:refined-semi-algebraic-2}] 
In the proof of Theorem \ref{thm:refined-semi-algebraic} instead of bounding
the number of semi-algebraically connected components of the various algebraic
sets 
\[
\tmop{Zer} (\{ Q_1, \ldots, Q_{\ell}, P_I \},
\R'^k)
\] 
using Theorem \ref{thm:refined-algebraic},
apply Theorem \ref{thm:refined-algebraic} directly to the sequence
$\mathcal{Q}, \tilde{\mathcal{P}}_I$, noting that its real zeros are the same
as $\tmop{Zer} (\{ Q_1, \ldots, Q_{\ell}, P_I \},
\R'^k)$, and also that the degree sequence
associated to $\tilde{\mathcal{P}}_I$ can be made to satisfy the requirement
of Theorem \ref{thm:refined-algebraic} by multiplying, for each $i$,  the $i$-th largest
degree in the sequence by $(k + 1)^{i - 1}$.
\end{proof}

\bibliographystyle{abbrv}
\bibliography{master}

\end{document}